\documentclass{article}
\usepackage{graphicx,amsmath,amsfonts,amssymb,amsthm,pdfpages,float} 
\usepackage{subcaption}
\usepackage{algorithm}
\usepackage{algorithmic}
\usepackage{comment}
\usepackage{todonotes}

\usepackage{authblk}

\usepackage[colorlinks=true, allcolors=blue]{hyperref}
\title{Low-rank approximation of Rippa method for RBF interpolation}
\author[1]{Jiawen Lyu\thanks{lyu.482@osu.edu}}
\author[1]{Maria Han Veiga\thanks{hanveiga.1@osu.edu}}
\affil[1]{Department of Mathematics, The Ohio State University}
\date{}
\usepackage[margin=1in]{geometry}

\newcommand{\R}{\mathbb{R}}

\begin{document}

\maketitle

\begin{abstract}
    We study the problem of selecting the shape parameter in Radial Basis function (RBF) interpolation  using leave-one-out-cross-validation (LOOCV). Since the classical LOOCV formula requires repeated solves with a dense $N \times N$ kernel matrix, we combine a Nystr\"{o}m approximation with the Woodbury identity to obtain an efficient surrogate objective that avoids large matrix inversions. Based on this reduced form, we compare a grid-based search with a gradient descent strategy and examine their behavior across different dimensions. Numerical experiments are performed in 1D, 2D, and 3D using the Inverse Multiquadratic RBF to illustrate the computational advantages of the approximation as well as the situations in which it may introduce additional sensitivity. These results show that the proposed acceleration makes LOOCV-based parameter tuning practical for larger datasets while preserving the qualitative behavior of the full method.
\end{abstract}

\section{Introduction}
Given a set of data points $\{x_{i}, f(x_{i})\}_{i = 1}^N $, where $x_{i}  \in \R^n$ and $f(x_{i}) \in \R$, the interpolation problem aims to find a function $s$ satisfying $s(x_{i}) = f(x_{i})$ for all $i=1,\cdots,n$.

Using radial basis functions (RBFs), we seek an approximation of the form
\begin{equation}
\label{eq:interpolation}
s(x) = \sum_{i = 1}^N \lambda_{i} \phi(\|x - x_{i}\|_{2}; \epsilon)
\end{equation}
where $\lambda_{i}\in\mathbb{R}$ are coefficients to be determined, $\phi$ is a positive definite radial basis function, and $\epsilon > 0$ denotes the \textit{shape parameter} that controls the shape of the basis.  Equivalently, \eqref{eq:interpolation} can be expressed as a system of equations
\begin{equation}
    \label{eq:interpolation_matrix}
    \mathbf{A}(X,\epsilon)\lambda = f(X),
\end{equation}
where $\mathbf{A}(x,\epsilon) \in \R^{N \times N}$ is the interpolation matrix with entries $\mathbf{A}_{ij} = \phi(\|x_{i}-x_{j}\|;\epsilon)$,  $X$ denotes the ordered set of interpolation nodes $\{x_{1}, \ldots, x_{N}\}$ and $f(X) = (f(x_{1}), \ldots, f(x_{N}))$ is the vector containing the unknown function evaluated at the interpolation nodes. 

The choice of $\epsilon$ plays a significant role in the quality of the interpolant $s$ \cite{mongillo}. For classical RBFs , small values of $\epsilon$ lead to ill-conditioned systems but high accuracy, while large values improve stability at the cost of precision \cite{schaback}.

Rippa \cite{rippa1999algorithm} introduced a widely used leave-one-out cross validation-based (LOOCV) method  that selects an optimal $\epsilon$ by minimizing the norm of an error vector $\vec{E}(\epsilon) = (E_{1}(\epsilon), \ldots, E_{N}(\epsilon))$:
\begin{equation}
\label{eq:rippa}
\arg \min_{\epsilon} \|(E_{1}(\epsilon), \ldots, E_{N}(\epsilon))\|, \text{ where } E_{k}(\epsilon) = f(x_{k}) - s^{(k)}(x_{k}, \epsilon) 
\end{equation}
with $s^{(k)}$ being the interpolant to a reduced data set, obtained by removing the node $x_{k}$ from $X$ and the corresponding function value $f(x_{k})$ from $f(X)$.

In order to solve the minimization \eqref{eq:rippa}, one often specifies a candidate search set for $\epsilon$ and 
the overall computational cost to evaluate the error $\vec{E}(\epsilon)$ for each candidate $\epsilon$ is $\mathcal{O}(N^3)$.

Building on the LOOCV method, Cavoretto et al. \cite{cavoretto2016rbf} proposed an RBF-PU scheme in which both the local shape parameter and the patch radius are chosen by minimizing a bivariate LOOCV error estimate on each subdomain, thereby improving robustness for non-uniform data distributions. Ghalichi et al. \cite{Ghalichi2022} modified the Rippa method by utilizing an auxiliary set of points to create a cost function that numerically mimics the behavior of the Root Mean Square *(RMS) error. Later Fasshauer and Zhang \cite{Fasshauer2007} extend it to apply in the setting of iterated approximate moving least squares approximation of
function value data and for RBF pseudo-spectral methods for the solution of
partial differential equations. The LOOCV method has also been used in several fields of science including optimization, machine learning, computational physics, and numerical analysis (e.g.\cite{Sarra2013,HOSSEINI2014602,DRISCOLL2002295}). More recently, Cavoretto et al. developed a Bayesian optimization approach for radial kernel parameter tuning \cite{CAVORETTO2024115716}, where the LOOCV-type error function is modeled by a Gaussian process and an acquisition strategy adaptively selects new candidate values of $\epsilon$, significantly reducing the number of expensive kernel solves required compared with exhaustive LOOCV grid search. 

Low-rank kernel approximation techniques, such as the Nystr\"{o}m method \cite{Williams2001} have been widely used to accelerate kernel-based algorithms by approximating the full kernel matrix  using a subset of representative landmark points. This idea has proven effective in large-scale kernel regression and Gaussian process inference \cite{Williams2001}, yet it has  not been fully explored in the context of RBF shape-parameter optimization.

In this paper, we address the two computational aspects of the classical LOOCV method:
\begin{itemize}
    \item \textbf{Naive evaluation of the error $\vec{E}(\epsilon)$ incurs a computational cost of $\mathcal{O}(N^3)$.} To address this bottleneck, we first derive an accelerated LOOCV objective function  using the Nystr\"{o}m approximation in conjunction with the Woodbury matrix identity \cite{Golub2013}. This new formulation avoids all $N \times N$ matrix operations, reducing the cost of a single LOOCV evaluation from $\mathcal{O}(N^3)$ to $\mathcal{O}(Nm^2 + m^3)$ and making the complexity linear with respect to $N$, enabling scalable shape parameter tuning for large datasets.
    \item \textbf{Traditional LOOCV requires the specification of a search range for $\epsilon$, which may exclude the optimal value.} To overcome this limitation, we apply a gradient-based optimization strategy that avoids the need to specify a search range.
\end{itemize}

We conduct extensive numerical experiments in one, two, and three dimensions to assess computational efficiency and approximation accuracy under various conditions. The remainder of this paper is organized as follows: Section \ref{sec:low-rank} presents the low-rank approximation; Section \ref{sec:gd} presents the gradient descent algorithm; Section \ref{sec:numerics} reports the numerical experiments and findings; and Section \ref{sec:conclusion} concludes with results and future directions.

\section{Low-Rank Approximation of the Interpolation Matrix}
\label{sec:low-rank}

For a fixed shape parameter $\epsilon$, evaluating the LOOCV objective requires solving \eqref{eq:interpolation_matrix}
where $\mathbf{A} \in \R^{N \times N}$ is the dense interpolation matrix built on the interpolation nodes set $X$, with entries \[\mathbf{A}_{ij} = \phi(\|x_{i}-x_{j}\|;\epsilon).\]
Throughout this work we consider the Inverse Multiquadric RBF
\[ \phi(r;\epsilon) = \frac{1}{\sqrt{1+(\epsilon r})^2},\]
although the techniques presented can be used for any positive definite RBFs. 

The cubic cost of forming and factorizing $\mathbf{A}$ makes repeated evaluations prohibitively expensive. To reduce this cost, we approximate $\mathbf{A}$ using a rank-$m$ Nystr\"{o}m expansion.

We observe that the interpolation matrix $\mathbf{A}$ depends on the interpolation nodes. We apply the k-means++ algorithm \cite{Arthur2007kmeans} to the interpolation nodes to identify a set of representative points $X_\ell$, which are given by the cluster centers, also called as \textit{landmarks} (see Appendix \ref{app:kmeans}). 

Then, given the $m$ landmarks selected by the k-means++ algorithm, we form the reduced matrices \[\mathbf{C} \in \R^{N \times m}, \quad  \mathbf{W} \in \R^{m \times m},\]
where  $\mathbf{C}_{ij} = \phi(\Vert x_i - x_j\Vert;\epsilon)$ for $x_i \in X$ and $x_j \in X_\ell$, and $\mathbf{W}_{ij} =\phi(\Vert x_i - x_j\Vert;\epsilon)$  for $x_i,x_j \in X_\ell$. The Nystr\"{o}m approximation is then \[\tilde{\mathbf{A}}_{nys}  = \mathbf{C} \mathbf{W}^{-1} \mathbf{C}^{\top},\]
where $W^{-1}$ denotes the inverse of $W$. Constructing the full interpolation matrix requires $\mathcal{O}(N^2)$ pairwise evaluations, followed by an $\mathcal{O}(N^3)$ factorization to compute LOOCV. Under Nystr\"{o}m approximation, the costs are reduced to $\mathcal{O}(Nm)$ and $\mathcal{O}(m^3)$, respectively.

\subsection{Woodbury-Accelerated LOOCV}
To evaluate LOOCV  efficiently, we first regularize the system \[\tilde{\mathbf{A}}_{reg} =\mathbf{C}\mathbf{W}^{-1} \mathbf{C}^{\top} + \lambda_{reg} I \]
for a small $\lambda_{reg}\in\mathbb{R}^+$ (we choose $\lambda_{reg} = 10^{-6}$ in all 1D, 2D, and 3D cases) and apply the Woodbury matrix identity \cite{Golub2013} to compute its inverse: \[\tilde{\mathbf{A}}_{reg}^{-1} = \lambda_{reg}^{-1}I  - \lambda_{reg}^{-1} C(W + \lambda_{reg}^{-1}C^{\top}C)^{-1} C^{\top} \lambda_{reg}^{-1}.\]
Doing so eliminates all $N \times N$  factorizations and reduces each LOOCV evaluation to $\mathcal{O}(Nm^2 +m^3)$, which becomes linear in $N$ once $m$ is fixed. 
\subsection{Landmark Stability and Random Seed Justification}
Since the k-means++ algorithm involves randomness in its initialization, one might  be concerned that fixing a single random seed  could inadvertently bias the numerical results. To verify that our conclusions do not depend on a particular initialization, we evaluate the stability of the landmark selection procedure by running the k-means++ algorithm for $R = 10$ times for each interpolation nodes dataset size $N$, with all runs using the same interpolation nodes but different cluster initial locations. We then compute the Normalized Mutual Information (NMI) \cite{Strehl2002ClusterEnsembles} between every pair of clustering results, measured as:
\[ NMI(C_1,C_2) = \frac{I(C_1,C_2)}{\sqrt{H(C_1)H(C_2)}},\]
where $C_1$ and $C_2$ are specific cluster assignments, $I(\cdot,\cdot)$ measures the mutual information. We compute the NMI between each pair of clustering results to quantify stability. For each configuration, we run k-means++  with $R = 10$ different random initialization seeds, yielding $\binom{10}{2} = 45$ pairwise comparisons. We report the mean and standard deviation of these NMI values. High mean NMI (close to 1) with low standard deviation indicates that the landmark selection is stable and insensitive to random initialization. Detailed computation steps are provided in Appendix~\ref{app:nmi}. 

In Figure \ref{fig:nmi_landmarks}, we report the mean and standard deviation of the NMI for dimensions $D = 1, 2, 3$ while varying the number of points $N$ and landmarks $m$. In all three cases, the NMI remains high, typically between 0.91 -- 0.94 in 1D, 0.86 -- 0.91 in 2D, and 0.82 -- 0.88 in 3D, even as $N$ increases from $1024$ to $8192$. Although a mild downward trend appears as the dimension and dataset size grow, the consistency across seeds remains strong: the standard deviations are small, and no outlier behavior is observed. These results confirm that the landmark selection is robust with respect to the randomness  in k-means++ initialization. In particular, the specific seed used in our main experiments is not a special or favorable -- the clustering structure it produces is representative of what k-means++ typically yields. This rules out the possibility of cherry-picking and ensures that the conclusions drawn from the Nystr\"{o}m-based LOOCV results are not sensitive to the choice of random seed. Furthermore, we note that the $NMI$ increases as the number of landmarks $m$ increases e.g., when $m = 50$, the NMI for 1D case stays around 0.9 and when $m = 400$, the NMI for 1D case is around 0.95.  

\begin{figure}
    \centering
    \begin{subfigure}{0.45\linewidth}
        \centering
        \includegraphics[width=\linewidth]{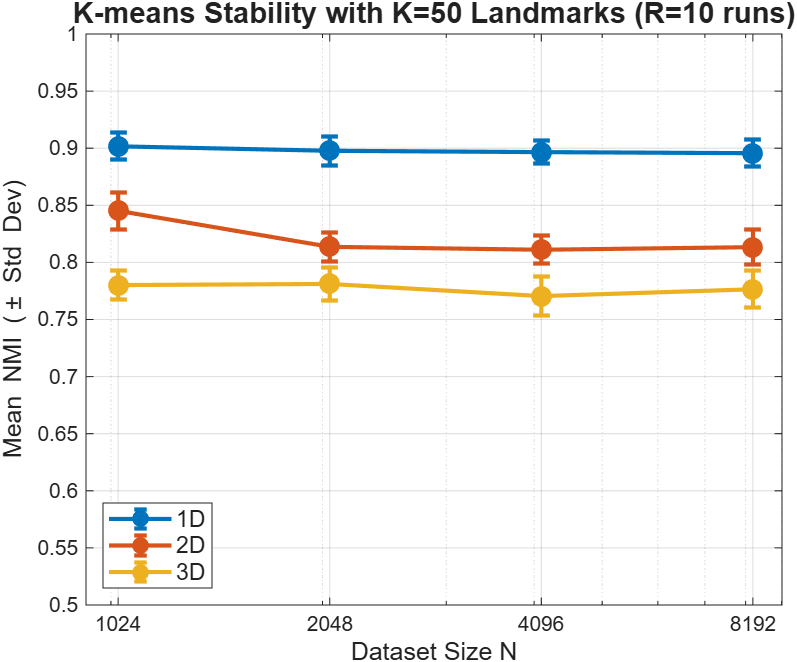}
        \caption{$m = 50$}
    \end{subfigure}
    \hfill
    \begin{subfigure}{0.45\linewidth}
        \centering
        \includegraphics[width=\linewidth]{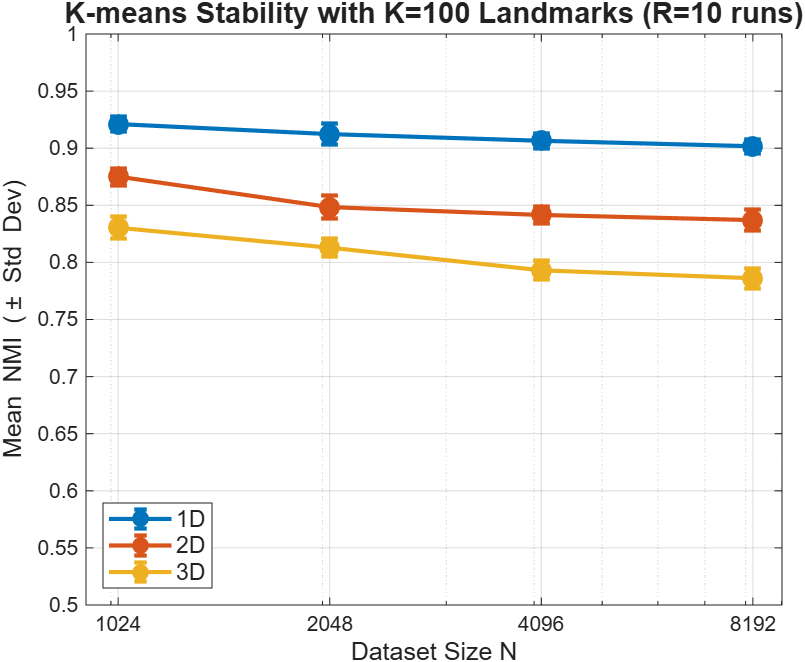}
        \caption{$m = 100$}
    \end{subfigure}
    \hfill
    \\
    \begin{subfigure}{0.45\linewidth}
        \centering
\includegraphics[width=\linewidth]{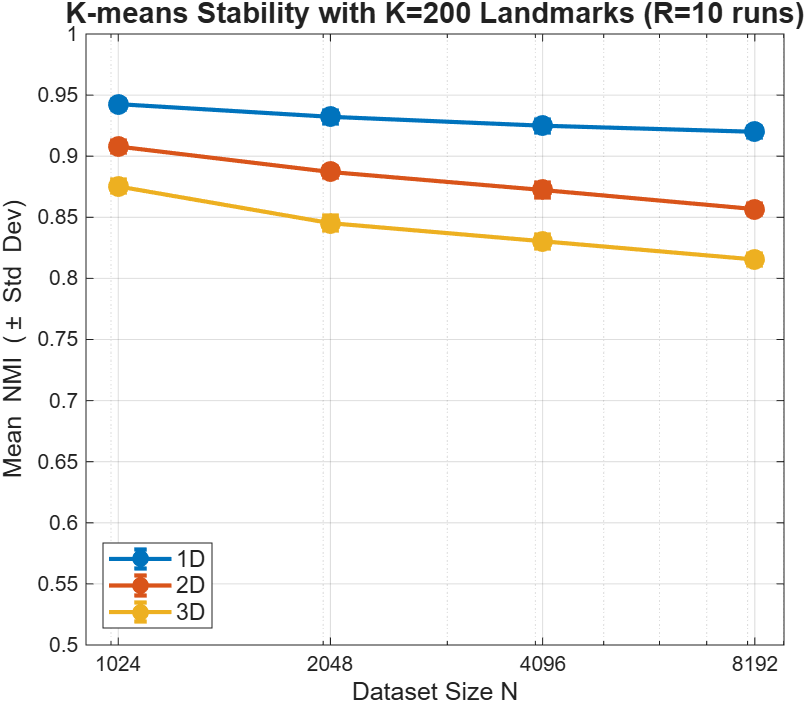}
        \caption{$m = 200$}
    \end{subfigure}
    \hfill
    \begin{subfigure}{0.45\linewidth}
        \centering
        \includegraphics[width=\linewidth]{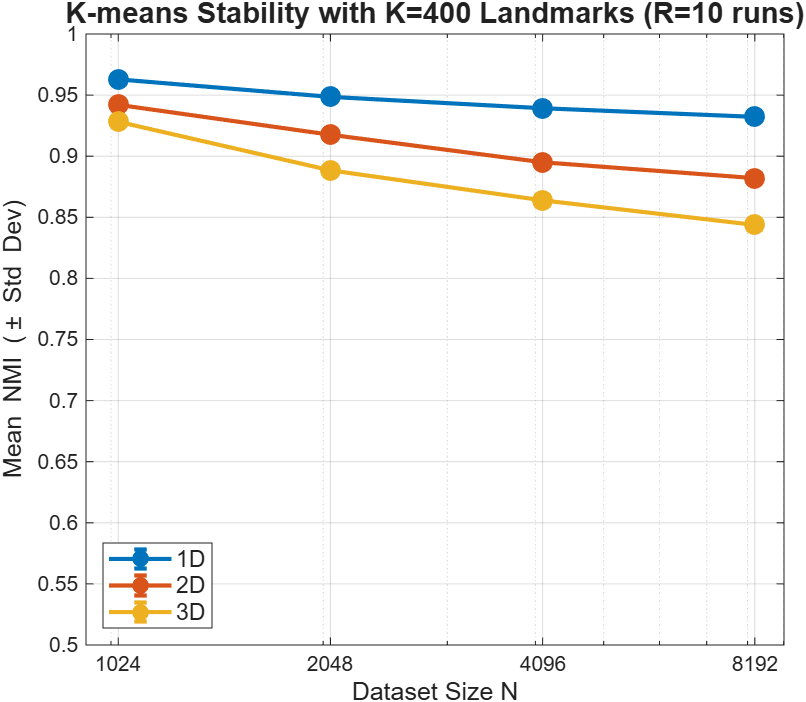}
        \caption{$m = 400$}
    \end{subfigure}
    \caption{NMI values for different numbers of landmarks}
    \label{fig:nmi_landmarks}
\end{figure}

\section{Gradient-Based Optimization of the Shape Parameter}
\label{sec:gd}
While Section \ref{sec:low-rank} addresses the computational cost of evaluating the LOOCV objective for a fixed shape parameter $\epsilon$, an efficient optimization strategy is still required to locate its minimizer. In this section, we investigate the use of a first-order optimization method -- gradient descent -- to solve the minimization \[\epsilon^* = \arg \min_{\epsilon > 0} L(\epsilon), \quad L(\epsilon):= \|E(\epsilon)\|_{2}^2,\]
where $E(\epsilon)$ is the vector of LOOCV residuals defined in \eqref{eq:rippa}. Since it is too costly to explicitly compute the analytical gradient of LOOCV objective with respect to $\epsilon$, numerical differentiation is employed.

To ensure the shape parameter remains strictly positive and to improve numerical stability across varying scales, we perform the optimization in the logarithmic space \cite{pedregosa2016hyperparameter}. Let $\theta = \ln (\epsilon)$, then, the optimization problem becomes finding $\theta^* = \arg \min_{\theta} L(e^{\theta})$. The derivative  with respect to the transformed variable is obtained via the chain rule: \[\frac{dL}{d\theta} = \frac{dL}{d\epsilon} \frac{d\epsilon}{d\theta} = \epsilon \cdot \frac{dL}{d\epsilon}.\]
Since the analytical derivative of the LOOCV objective is computationally prohibitive, we approximate $dL/d\epsilon$ using a centered finite-difference scheme: \[\frac{dL}{d\epsilon} \approx \frac{L(\epsilon + h) - L(\epsilon - h)}{2h},\]
where the step size is chosen  adaptively as $h = \max \{10^{-8},10^{-8}|\epsilon|\}$ to balance numerical precision and stability.

\subsection{Adaptive Step Size and Line Search}
Since the landscape of the LOOCV objective function can be highly sensitive, characterized by narrow valleys and regions that trigger numerical overflow, a standard fixed-step gradient descent is prone to divergence. To address this, we employ a robust backtracking line search strategy combined with a heuristic to dynamically adapt the learning rate $\eta$. At each iteration $k$, given the current log-parameter $\theta_k = \ln(\epsilon_k)$ and the gradient $\nabla_\theta L$, we propose a candidate update:$$\theta_{candidate} = \theta_k - \eta_k \nabla_\theta L$$
We then evaluate the objective function at the corresponding shape parameter $\epsilon_{candidate} = \exp(\theta_{candidate})$. The step is accepted or rejected based on the following criteria: first, we check if the candidate loss $L(\epsilon_{candidate})$ yields a valid numerical result. Due to the ill-conditioning of RBF matrices, aggressive steps can occasionally result in singular matrices, leading to NaN or Inf values in the objective. Such steps are immediately treated as failures and the step is rejected. If the value is valid, we check for a strict decrease in the objective function, i.e., $L(\epsilon_{candidate}) < L(\epsilon_k)$.
Based on the outcome, the step size $\eta$ is updated adaptively \cite{vogl1988accelerating}: we set the initial step size $\eta$ to be $1.0$.  If the step is accepted, we update $\theta_{k+1} = \theta_{candidate}$. To accelerate convergence in flat regions of the landscape, we slightly increase the step size for the next iteration: $\eta_{k+1} = \min(1.2 \eta_k, \eta_{max})$, where $\eta_{max}=5.0$ prevents excessive growth. If the objective increases or is invalid, the step is rejected. We perform a backtracking procedure within the current iteration by reducing the step size: $\eta_k = 0.5 \eta_k$. This reduction is repeated up to a maximum of 15 retries. If no valid descent direction is found after these attempts, the algorithm terminates to avoid stagnation. 

This adaptive scheme ensures that the optimization cautiously navigates the ill-conditioned regions while moving rapidly when the gradient information is reliable.

\subsection{Initialization and Stopping Criterion}
Since the candidate set for $\epsilon$ has been removed, the gradient descent algorithm requires an explicit initial guess. To obtain a reliable starting point for optimization, for 1D case in the numerical experiments, we adopt the initialization strategy introduced by Richard Franke \cite{Franke1982}, where the initial shape parameter is chosen as \[\epsilon = 0.8\sqrt{N}/D,\]
where $D$ is the diameter of the minimal circle enclosing all data points and $N$ is the size of interpolation nodes. In 2D case, we adopt the modified Franke method \cite{kuo2015selection} to deal with the numerical instability because of increase of dimension, where the initial shape parameter is chosen as \[\epsilon = 0.8N^{1/4}/D.\] In 3D, this rule became less robust due to the stronger ill-conditioning of the kernel matrix. Instead, we use a data-driven distance heuristic \cite{Jaakkola1998Exploiting}: from at most $1000$ randomly selected training points, we compute all pairwise Euclidean distances, take their median $d_{med}$, and set \[\epsilon_{0} = \frac{1}{d_{med} +10^{-8}}.\]
This choice ties the initial shape parameter to the typical node spacing and has proved stable across all 3D test cases.
We terminate the optimization when the magnitude of $\frac{dL}{d\epsilon}$  falls below $10^{-8}$, or when the relative change in $L(\epsilon)$ becomes smaller than $10^{-12}$ over consecutive iterations. In all experiments, we also enforce a maximum of 100 iterations.

\section{Numerical Experiments}
\label{sec:numerics}
In this section, we present numerical experiments aimed at comparing different strategies for selecting the shape parameter in radial basis function (RBF) interpolation. Specifically, we consider the classical Rippa method (Rippa), the reformulated Rippa method accelerated using the Nystr\"{o}m approximation (Rippa + Nystr\"{o}m), the standard gradient descent approach (GD), and gradient descent with Nystr\"{o}m acceleration (GD + Nystr\"{o}m). We present experiments in one, two, and three spatial dimensions. The interpolation nodes are uniformly sampled from $[0,1]^d$, with $d = 1,2,3$. All methods are tested under identical experimental conditions to ensure a fair comparison.

The evaluation of the methods' errors are done in the following way: let $f$ be a target function defined on a $d$-dimensional domain. Given a set of $N$ interpolation nodes
\[
X = \{(x_{i1}, \ldots, x_{id})\}_{i=1}^{N} \subset [0,1]^d,
\]
we define the data vector $f(X) \in \mathbb{R}^N$ by evaluating $f$ at the interpolation nodes.

For a fixed value of the shape parameter $\epsilon$, we solve the RBF interpolation problem \eqref{eq:interpolation} and obtain the corresponding interpolant $s$. The accuracy of the approximation is assessed using the root-mean-square (RMS) error,
\[
\mathrm{RMS}(\epsilon)
= \sqrt{\frac{1}{N_{\text{test}}}
\sum_{i=1}^{N_{\text{test}}}
\bigl[s(x_{i1}, \ldots, x_{id}) - f(x_{i1}, \ldots, x_{id})\bigr]^2}.
\]
The test points are uniformly sampled over $[0,1]^d$ and fixed across all experiments for a given test function.  Unless otherwise stated, the number of test points is fixed to $N_{\text{test}} = 5000$. The optimal shape parameter $\epsilon^\ast$ is defined as the value of $\epsilon$ that minimizes the RMS error.

\paragraph{Shape Parameter Search Strategy.} For the Rippa based methods,  we employ a unified two-stage coarse-to-fine search strategy across all experiments to identify the optimal shape parameter. In the first stage, a coarse search is conducted over $30$ logarithmically spaced values in the interval $[10^{-5},10^{3}]$, yielding an approximate minimizer $\epsilon_{\text{coarse}}$. In the second stage, the search is refined over the interval $[0.5\,\epsilon_{\text{coarse}},\,2.0\,\epsilon_{\text{coarse}}]$ using $50$ linearly spaced points. This adaptive strategy provides a good balance between robustness, accuracy, and computational efficiency.

\paragraph{Interpolation Sizes and Nystr\"{o}m Parameters.} For the standard (non-Nystr\"{o}m) methods, the number of interpolation nodes is varied over
\[
N \in \{64, 128, 256, 512, 1024, 2048, 4096\}.
\]
For the Nystr\"{o}m-accelerated methods, we consider
\[
N \in \{512, 1024, 2048, 4096\}.
\]
In all Nystr\"{o}m experiments, the number of landmark points is fixed to $m = 200$. From Figure \ref{fig:nmi_landmarks}, we could see that as $m$ increases from $50$ to $200$, the NMI in all three dimensions increases by $0.15$. However, as $m$ increases from $200$ to $400$, the NMI in all three dimensions increases by only $0.005$, which is less noticeable. So as a  trade-off between accuracy and computational efficiency, we fix the number of landmarks to be $m = 200$ throughout our numerical experiments. 

\paragraph{Computational Details and Regularization.}
All numerical experiments are carried out in \textsc{MATLAB} on an Intel(R) Core(TM) i5-12500 CPU running at 2.50\,GHz. To improve numerical stability, the interpolation systems are regularized by adding a small multiple of the identity matrix. Specifically, we use a regularization parameter $\lambda_{\text{reg}} = 10^{-14}$ in the one-dimensional case and $\lambda_{\text{reg}} = 10^{-10}$ in the two- and three-dimensional cases.
Code is available on \url{https://github.com/Christine523/low-rank-approximation-RBF}.

\subsection{One-dimensional case}

We first evaluate the performance of all four methods on two one-dimensional test functions: the smooth, transcendental function $f_1(x)$ and the classic Runge function $f_2(x)$.\\

The test function $f_{1}: [0,1] \rightarrow \R$ is \[f_{1}(x) = e^{\sin(\pi x)}.\]

\begin{figure}[H]
    \centering
    \begin{subfigure}[b]{0.5\textwidth} 
        \centering
        \includegraphics[width=\linewidth]{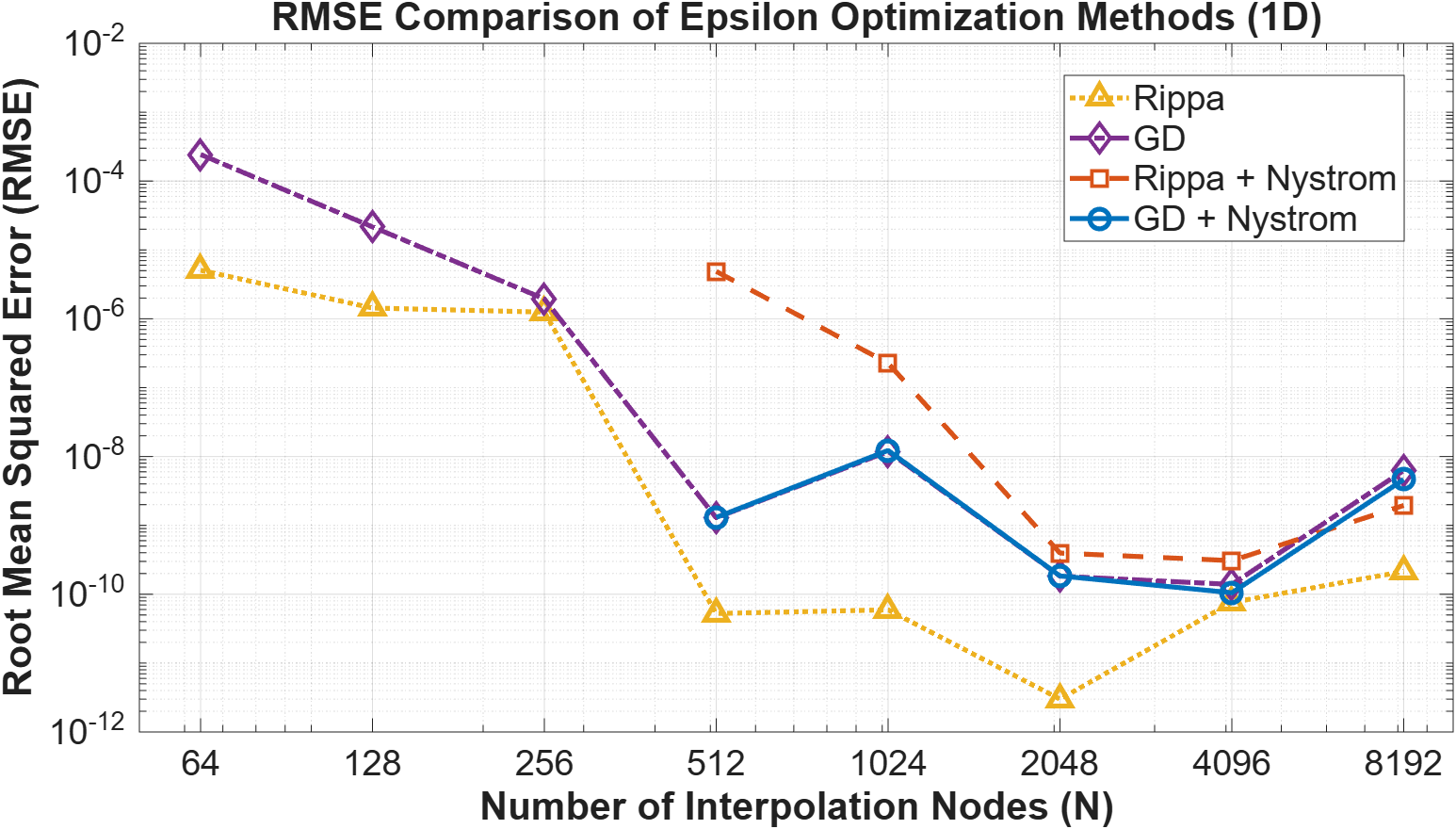}
        \caption{Accuracy Comparison}
        \label{fig:f1a}
    \end{subfigure}
    
    \vspace{1em} 

    \begin{subfigure}[b]{0.5\textwidth}
        \centering
        \includegraphics[width=\linewidth]{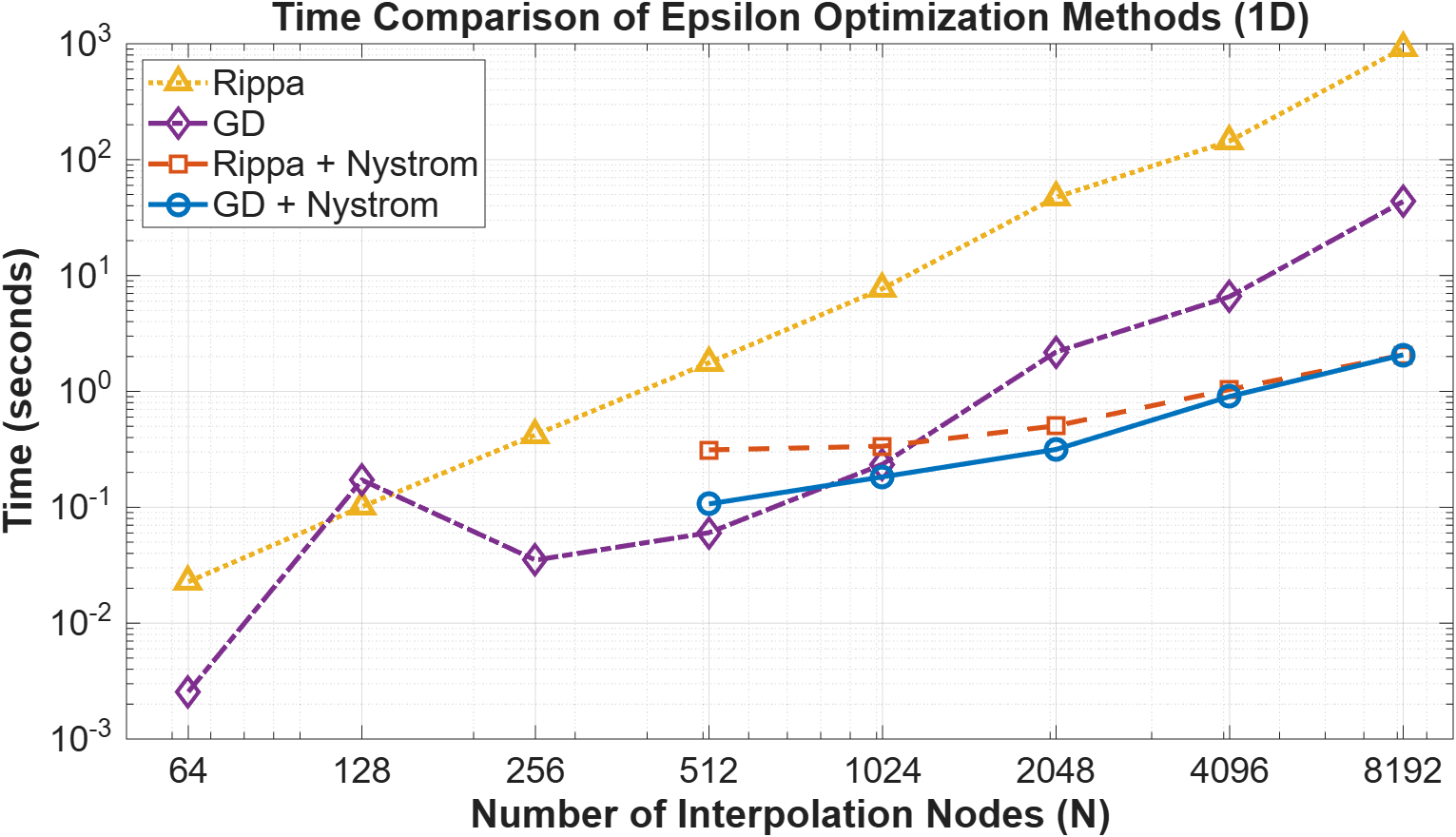}
        \caption{Time comparison}
        \label{fig:f1b}
    \end{subfigure}

    \caption{Performance plot for function $f_{1}$}
    \label{fig:f1plot}
\end{figure}

The test function $f_{2}: [0,1]\rightarrow \R$ is \[f_{2}(x) = \frac{1}{1+16x^2}.\]

\begin{figure}[H]
    \centering
    \begin{subfigure}[b]{0.7\textwidth} 
        \centering
        \includegraphics[width=\linewidth]{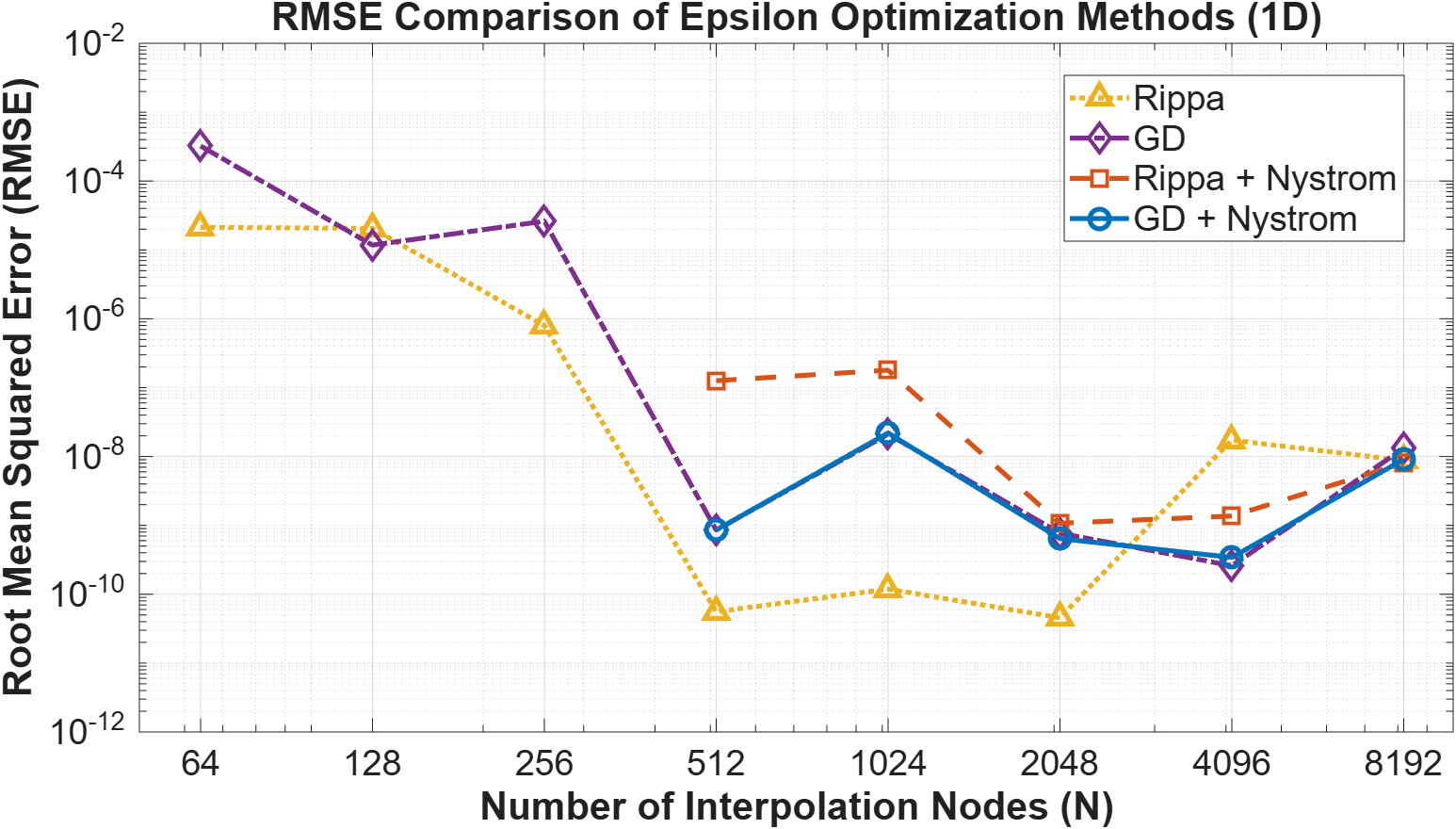}
        \caption{Accuracy Comparison}
        \label{fig:f2a}
    \end{subfigure}
    
    \vspace{1em} 

    \begin{subfigure}[b]{0.7\textwidth}
        \centering
        \includegraphics[width=\linewidth]{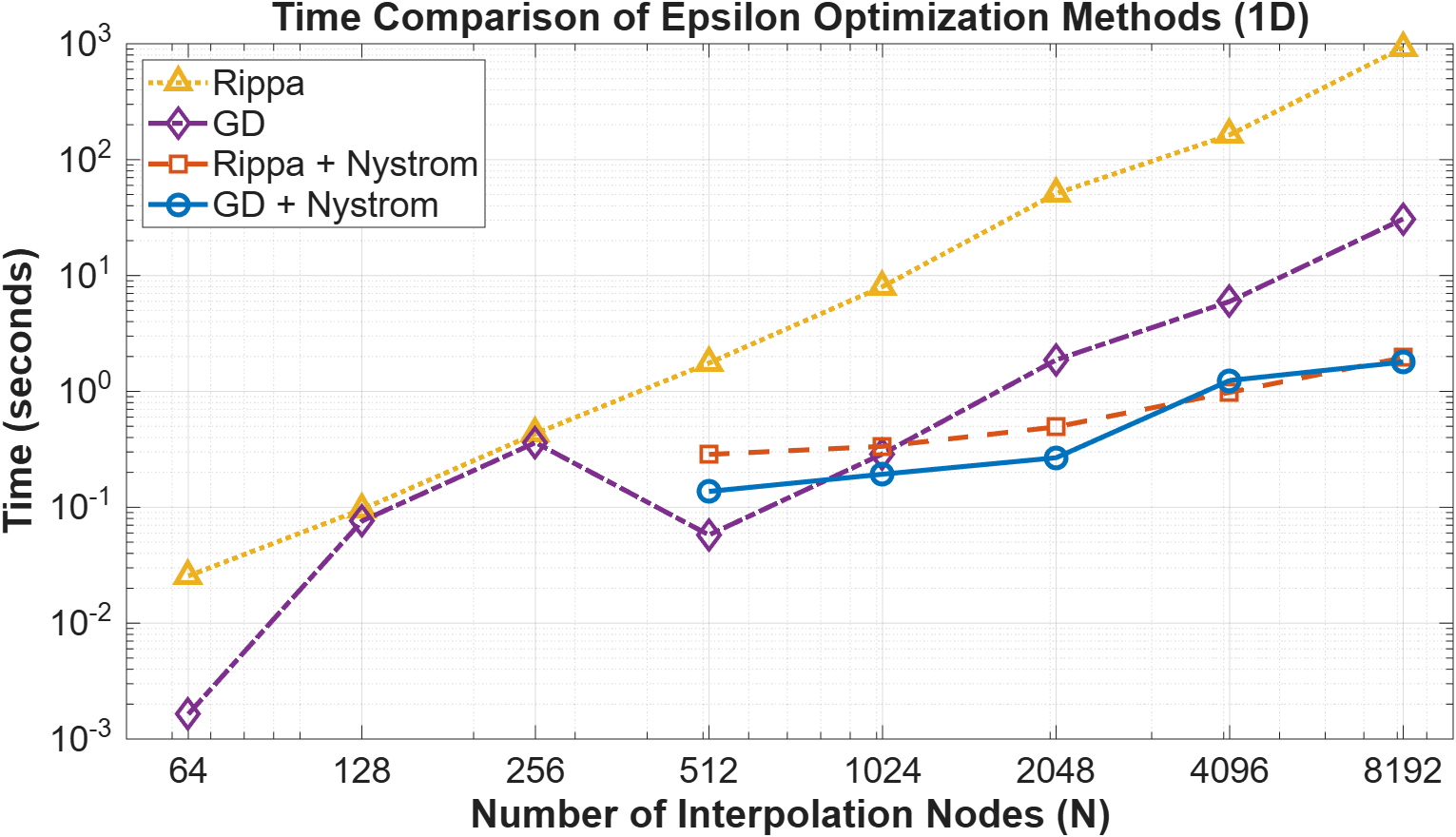}
        \caption{Time comparison}
        \label{fig:f2b}
    \end{subfigure}

    \caption{Performance plot for function $f_{2}$}
    \label{fig:f2plot}
\end{figure}

The timing results for these 1-dimensional test cases (see Figures \ref{fig:f1b} and \ref{fig:f2b}) demonstrate the efficiency of the Nystr\"{o}m formulation, with both accelerated methods scaling linearly. Notably, ``GD with Nystr\"{o}m" is the fastest method, outperforming the grid-based ``Rippa with Nystr\"{o}m" by requiring significantly fewer objective evaluations to converge. However, the accuracy of these methods (see Figures \ref{fig:f1a} and \ref{fig:f2b}) reflect the inherent instability of 1D RBF interpolation. The severe ill-conditioning of the interpolation matrix creates a sensitive optimization landscape, causing all methods to exhibit occasional error spikes. Here, ``Gradient Descent with Nystr\"{o}m" method slightly outperforms ``Rippa with Nystr\"{o}m" method. 

\subsection{Two-dimensional case}
\begin{figure}[H]
    \centering
    \begin{subfigure}[b]{0.7\textwidth} 
        \centering
        \includegraphics[width=\linewidth]{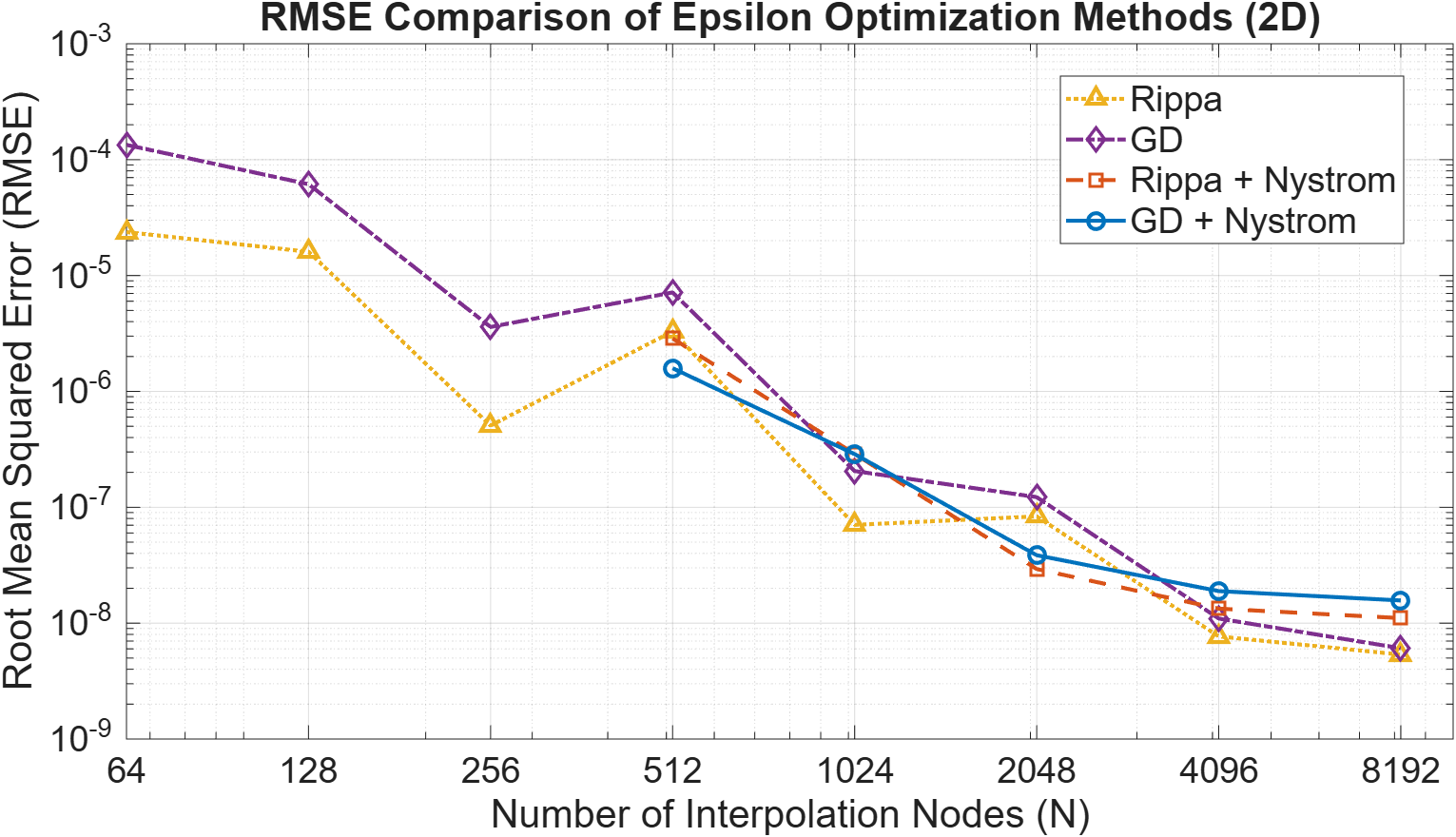}
        \caption{Accuracy Comparison}
        \label{fig:f3a}
    \end{subfigure}
    
    \vspace{1em} 

    \begin{subfigure}[b]{0.7\textwidth}
        \centering
        \includegraphics[width=\linewidth]{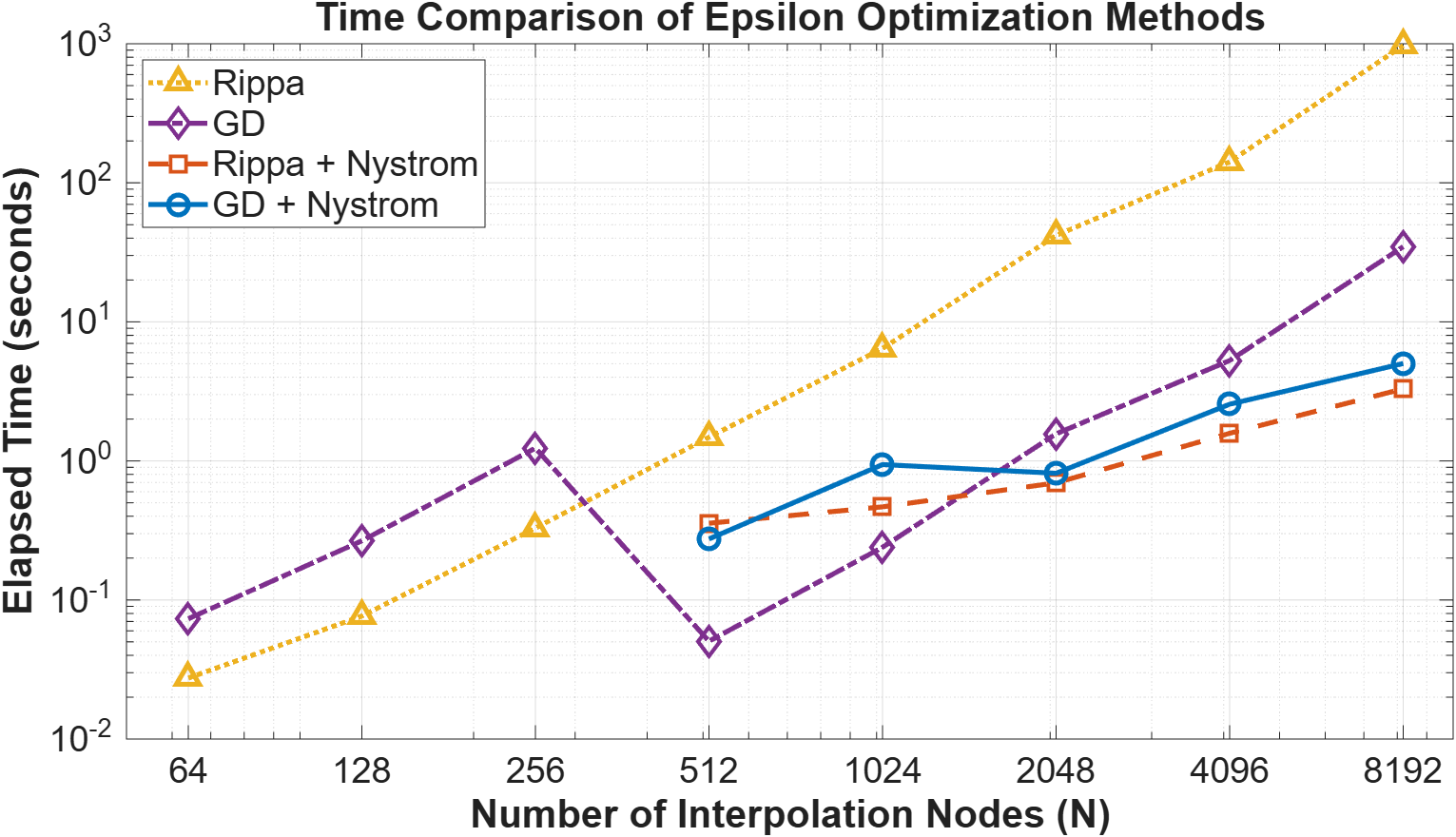}
        \caption{Time comparison}
        \label{fig:f3b}
    \end{subfigure}

    \caption{Performance plot for function $f_{3}$}
    \label{fig:f3plot}
\end{figure}

\begin{figure}[H]
    \centering
    \begin{subfigure}[b]{0.7\textwidth} 
        \centering
        \includegraphics[width=\linewidth]{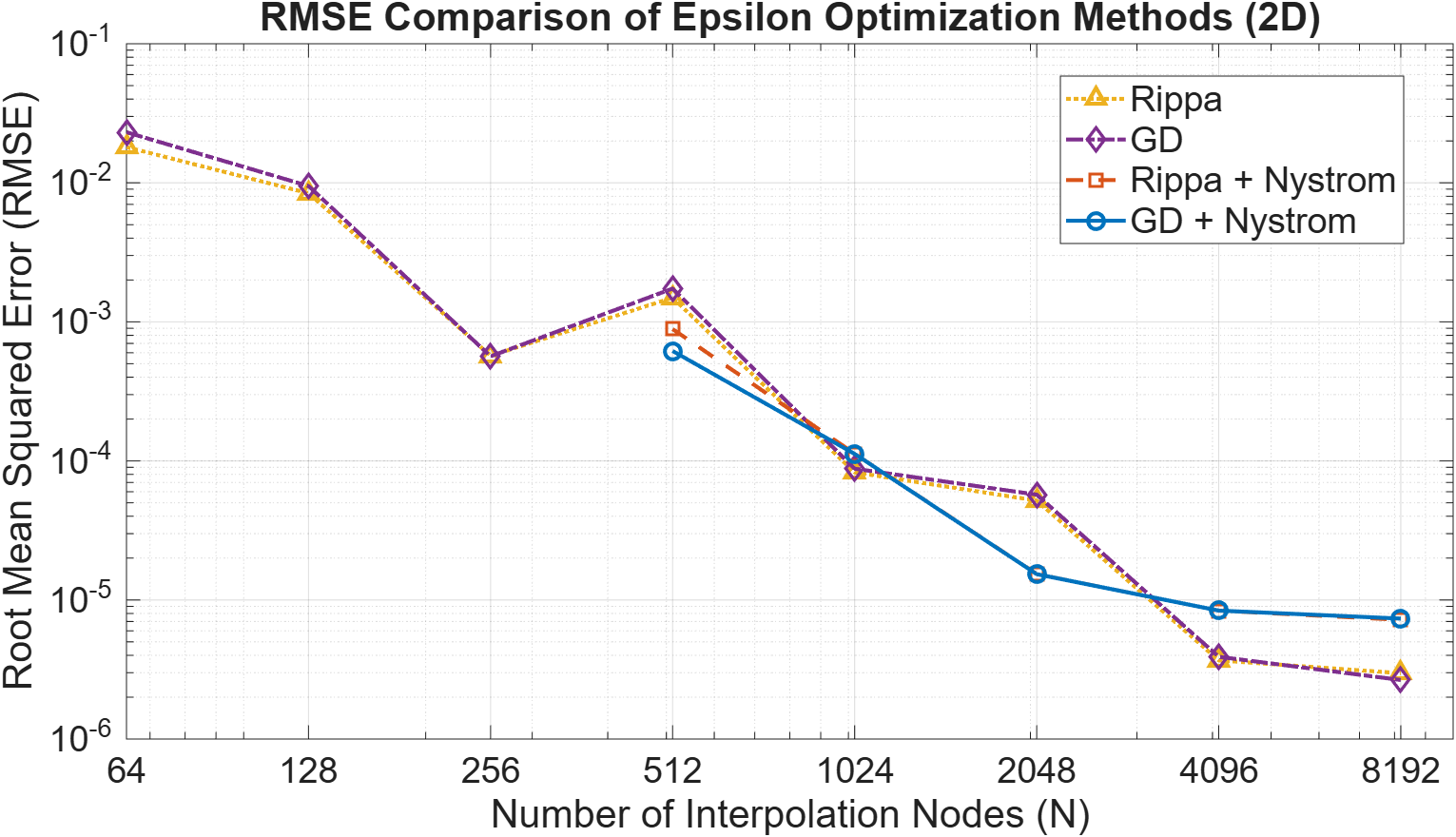}
        \caption{Accuracy Comparison}
        \label{fig:f4a}
    \end{subfigure}
    
    \vspace{1em} 

    \begin{subfigure}[b]{0.7\textwidth}
        \centering
        \includegraphics[width=\linewidth]{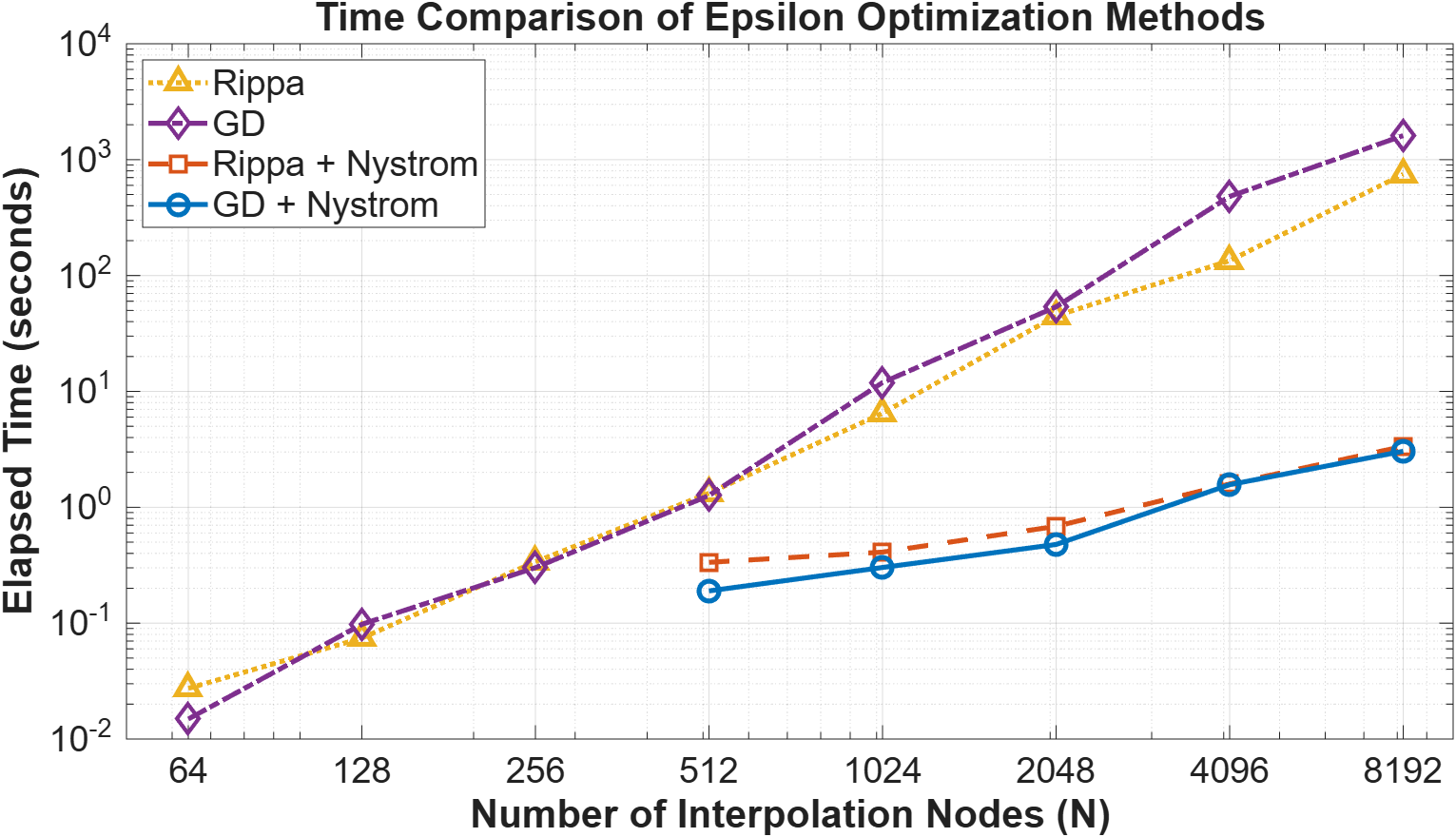}
        \caption{Time comparison}
        \label{fig:f4b}
    \end{subfigure}

    \caption{Performance plot for function $f_{4}$}
    \label{fig:f4plot}
\end{figure}

\begin{figure}[H]
    \centering
    \begin{subfigure}[b]{0.7\textwidth} 
        \centering
        \includegraphics[width=\linewidth]{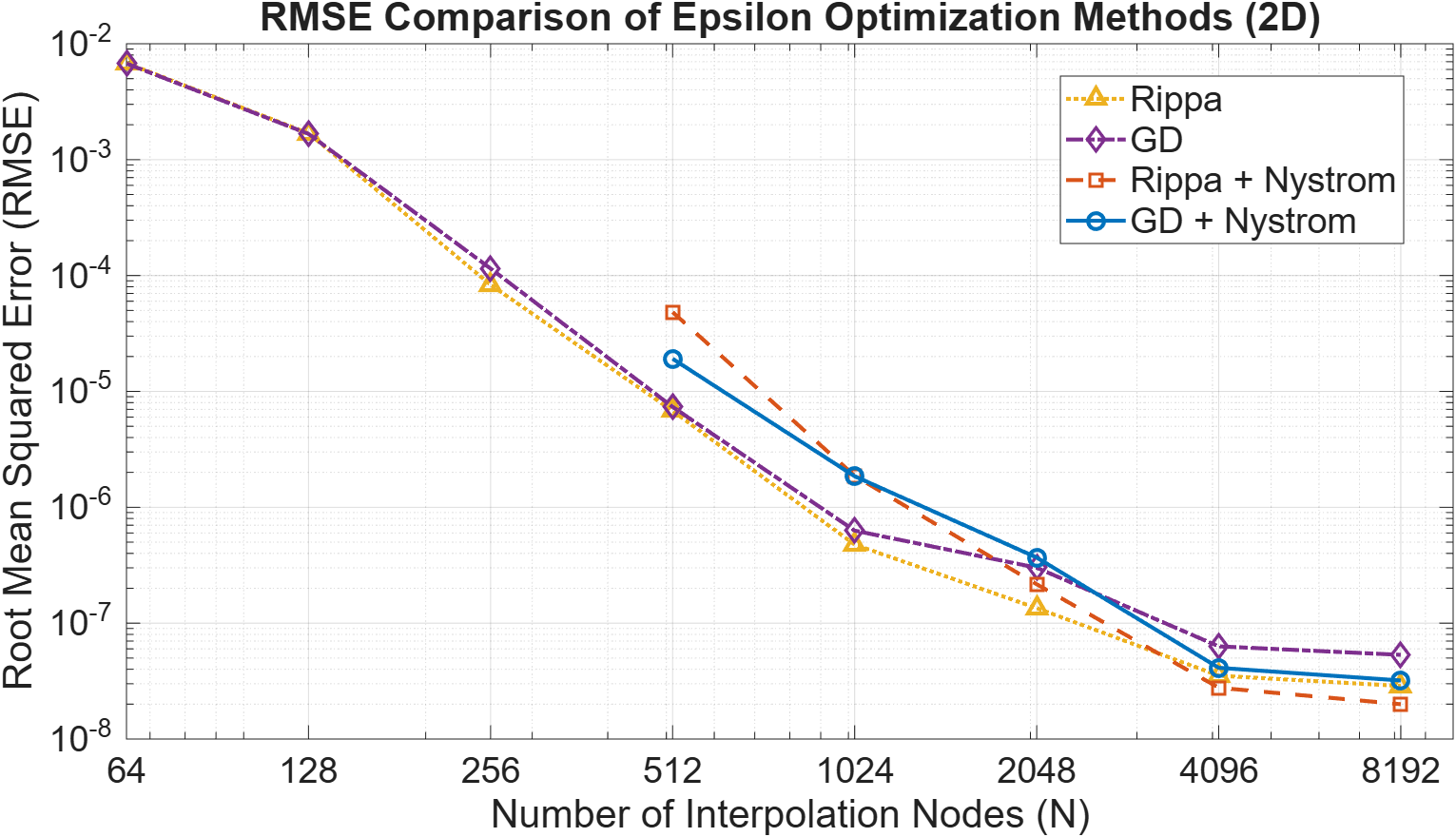}
        \caption{Accuracy Comparison}
        \label{fig:f5a}
    \end{subfigure}
    
    \vspace{1em} 

    \begin{subfigure}[b]{0.7\textwidth}
        \centering
        \includegraphics[width=\linewidth]{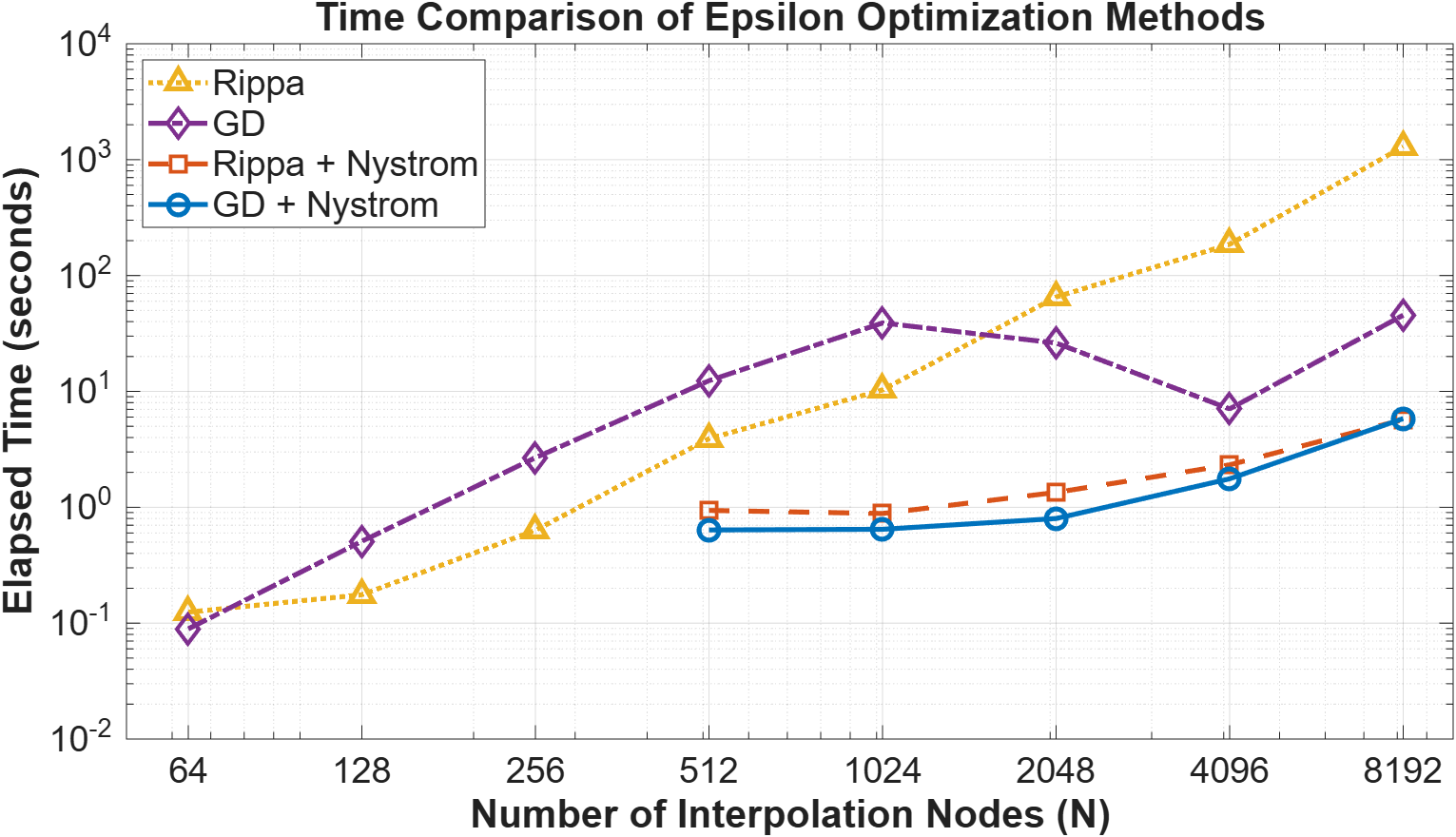}
        \caption{Time comparison}
        \label{fig:f5b}
    \end{subfigure}

    \caption{Performance plot for function $f_{5}$}
    \label{fig:f5plot}
\end{figure}

We further validate our findings using two-dimensional test functions: $f_{3}(\mathbf{x})$ (a smooth function), $f_{4}(\mathbf{x})$ (a steeper smooth function), and $f_{5}(\mathbf{x})$ (the famous Franke's function). The performance plots are shown in Figures \ref{fig:f3plot}, \ref{fig:f4plot} and \ref{fig:f5plot} respectively.

The test function $f_{3}: [0,1]^2 \rightarrow \R$ is \[f_{3}(\mathbf{x}) = (1+e^{-1}-e^{-x_{1}} - e^{-(x_{1}-1)})(1 +e^{-1}-e^{-x_{2}}-e^{-(x_{2}-1)}).\]

The test function $f_{4}:  [0,1]^2  \rightarrow \R$ is

\[f_{4}(\mathbf{x}) = (1+e^{-\frac{1}{0.1}} - e^{-\frac{x_{1}}{0.1}} - e^{\frac{(x_{1}-1)}{0.1}})(1 +e^{-\frac{1}{0.1}} - e^{-\frac{x_{2}}{0.1}} - e^{\frac{(x_{2}-1)}{0.1}}).\]

The test function $f_{5}: [0,1]^2  \rightarrow \R$ is the Franke's function: 

\[f_{5}(\mathbf{x}) = \frac{3}{4}e^{-(\frac{(9x_{1}-2)^2 +(9x_{2}-2)^2}{4})} +\frac{3}{4} e^{-(\frac{(9x+1)^2}{49} + \frac{(9x_{2} + 1)^2}{10})} + \frac{1}{2} e^{-(\frac{(9x-7)^2 +(9y-2)^2}{4})} - \frac{1}{5}e^{-((9x-4)^2) + (9y-7)^2}.\]

The two-dimensional experiments demonstrate significantly improved stability compared to the 1-dimensional case. The timing results (see Figures \ref{fig:f3b}, \ref{fig:f4b}, and \ref{fig:f5b}) confirm the linear scalability of the Nystr\"{o}m-based approaches. Consistent with the 1-dimensional findings, ``GD with Nystr\"{o}m" achieves the lowest runtime, surpassing the grid-based ``Rippa with Nystr\"{o}m" by efficiently navigating the optimization landscape with fewer objective evaluations.

In this case, the ``Gradient Descent with Nystr\"{o}m" method achieves lower RMSE than standard LOOCV when $N = 2048, 4096, 8192$. 

\subsection{Three-dimensional case}

To assess scalability  in higher dimensions, we evaluate the four methods on three 3D test functions, $f_{6}, f_{7}$, and $f_{8}$ \cite{Fornberg2011Stable}, defined on the unit cube domain $[0,1]^3$. As the dimension increases, kernel matrices typically become more ill-conditioned and interpolation becomes more challenging. This setting provides an important stress-test for the proposed acceleration framework.

Here, we modify  the training node distribution to ensure numerical stability. While the 1D and 2D tests use uniformly sampled points, such nodes are known to produce severely ill-conditioned interpolation matrices in higher dimensions. This behavior is consistent with the result of Platte, Trefethen, and Kuijlaars \cite{Platte2011Impossibility}, who showed that approximation schemes based on equispaced samples cannot remain stable as the dimension or resolution increases. 

To avoid this instability, the training points in 3D are transformed using a cosine mapping that concentrates more nodes near the boundary \[x_k  \mapsto 0.5 \times (1-\cos(\pi x_k)), \hspace{0.5cm} k \in \{1,2,3\} \text{ and } \mathbf{x} = (x_1, x_{2},x_{3}),\]

mimicking Chebyshev-type clustering. This adjustment improves the conditioning of the interpolation matrix and allows Nystr\"{o}m-Woodbury LOOCV formulation to produce meaningful and comparable results. The test points remain uniformly distributed, ensuring that the accuracy evaluation is not biased by the choice of training nodes.

The test function $f_{6}: \Omega_{3} \rightarrow \R$ is \[f_{6}(\mathbf{x}) = 1\]

\begin{figure}[H]
    \centering
    \begin{subfigure}[b]{0.7\textwidth} 
        \centering
        \includegraphics[width=\linewidth]{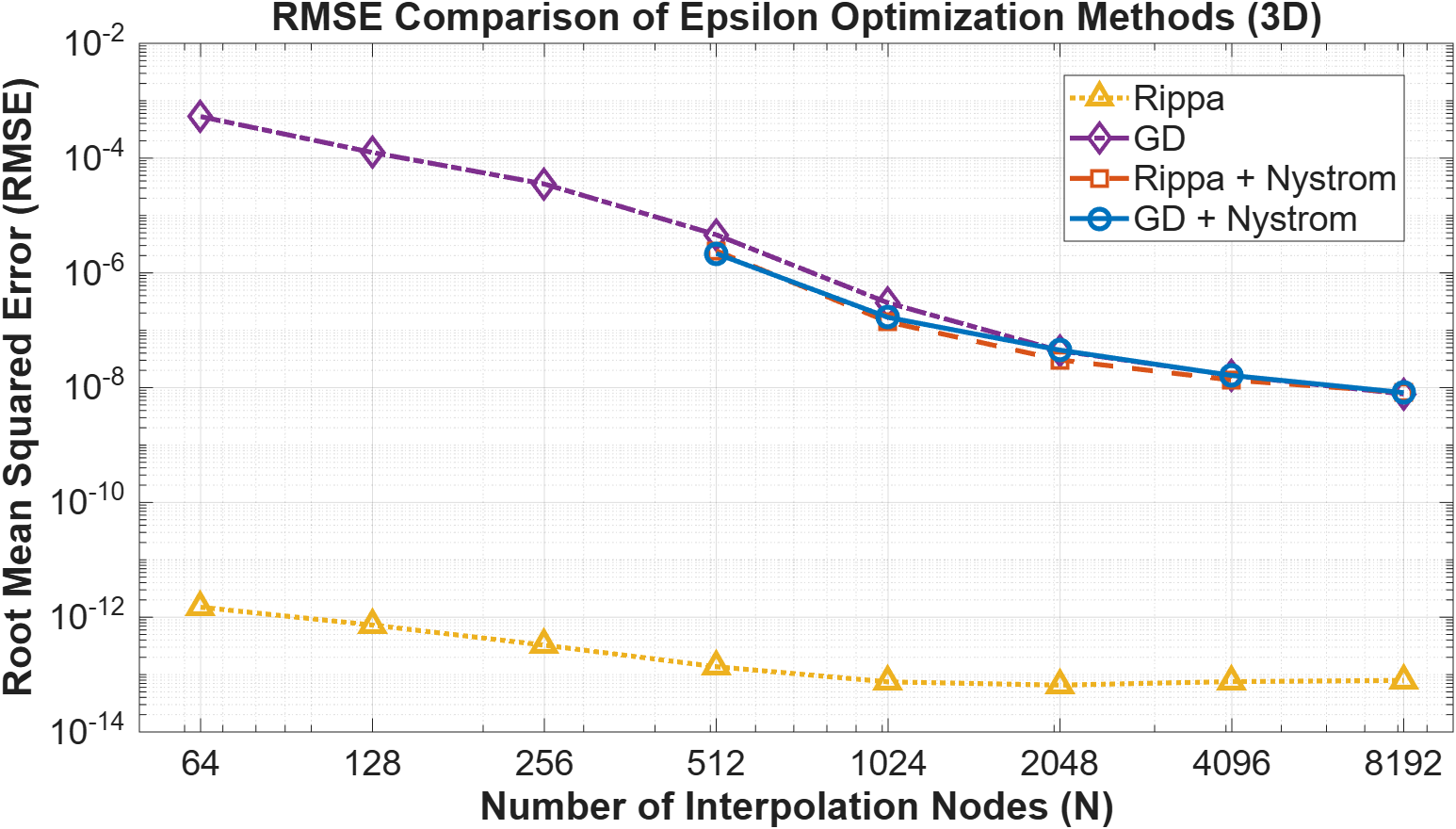}
        \caption{Accuracy Comparison}
        \label{fig:f6a}
    \end{subfigure}
    
    \vspace{1em} 

    \begin{subfigure}[b]{0.7\textwidth}
        \centering
        \includegraphics[width=\linewidth]{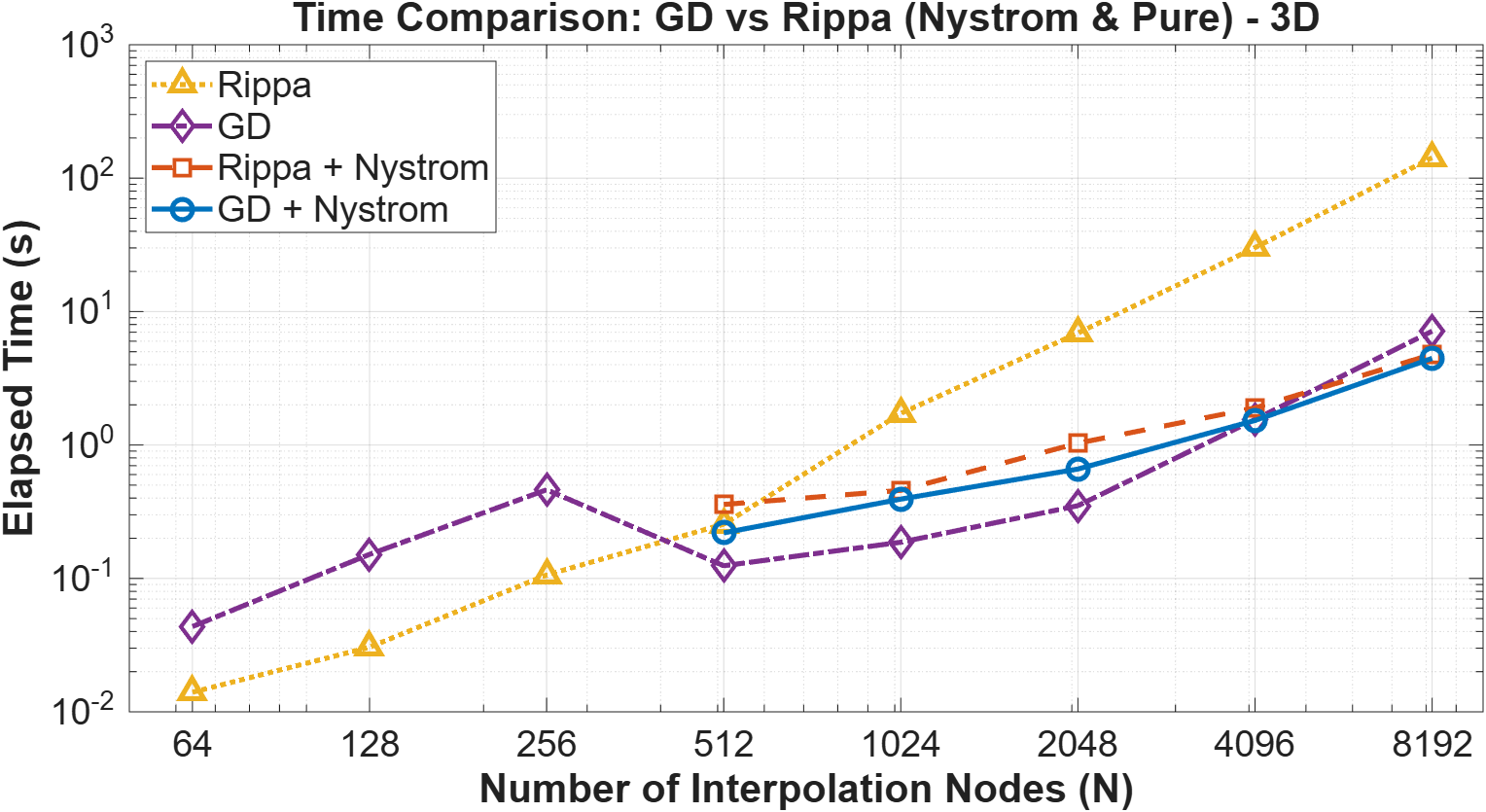}
        \caption{Time comparison}
        \label{fig:f6b}
    \end{subfigure}

    \caption{Performance plot for function $f_{6}$}
    \label{fig:f6plot}
\end{figure}

The test function $f_{7}: \Omega_{3} \rightarrow \R$ is \[f_{7}(\mathbf{x}) = \sin(x_{1}^2 + 2x_{2}^2) - \sin(2x_{1}^2 + (x_{2}-0.5)^2 + x_{3}^2),\]
$(x_{1},x_{2},x_{3}) \in [0,1] \times [0,1] \times [0,1]$.

\begin{figure}[H]
    \centering
    \begin{subfigure}[b]{0.7\textwidth} 
        \centering
        \includegraphics[width=\linewidth]{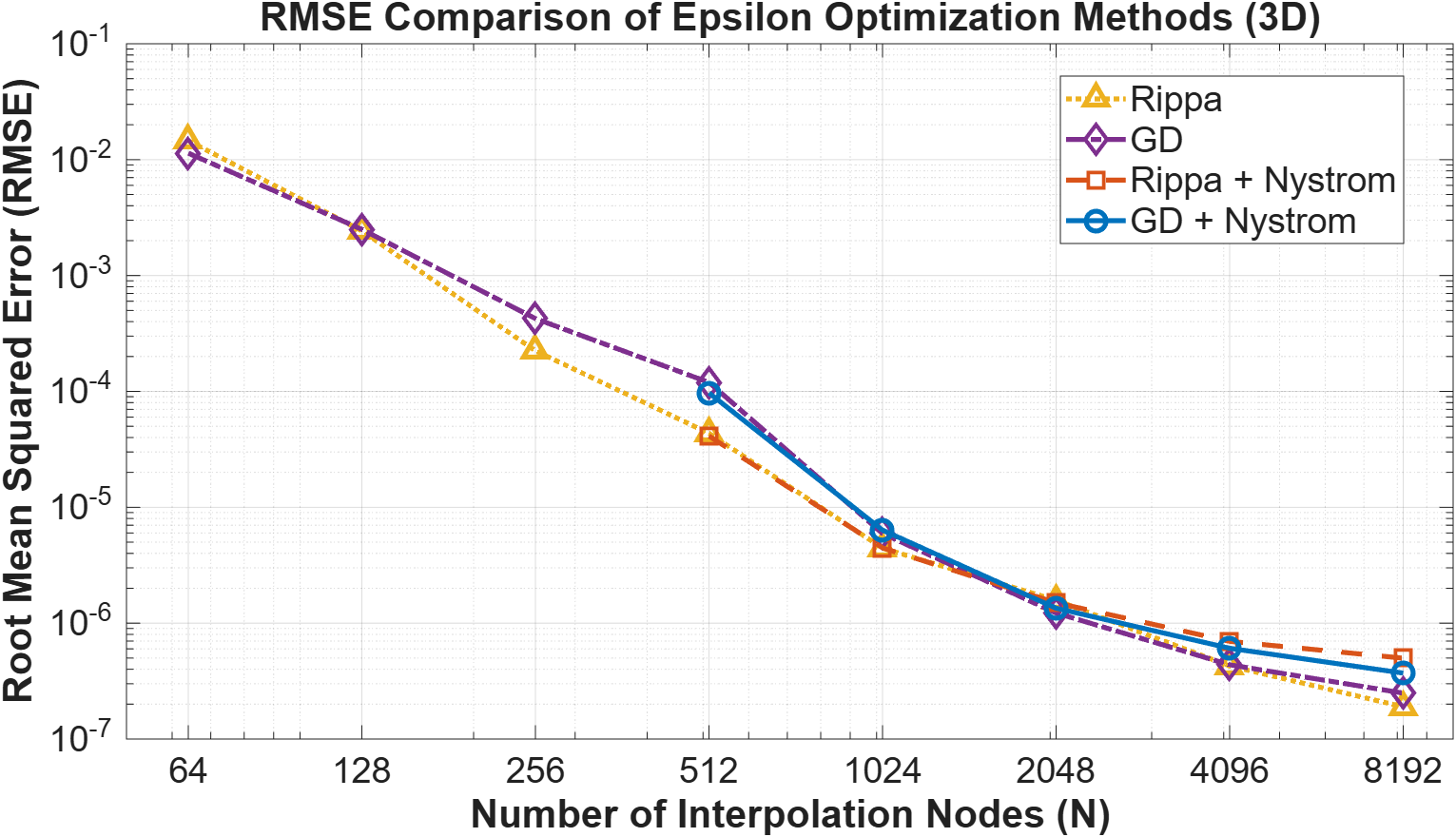}
        \caption{Accuracy Comparison}
        \label{fig:f7a}
    \end{subfigure}
    
    \vspace{1em} 

    \begin{subfigure}[b]{0.7\textwidth}
        \centering
        \includegraphics[width=\linewidth]{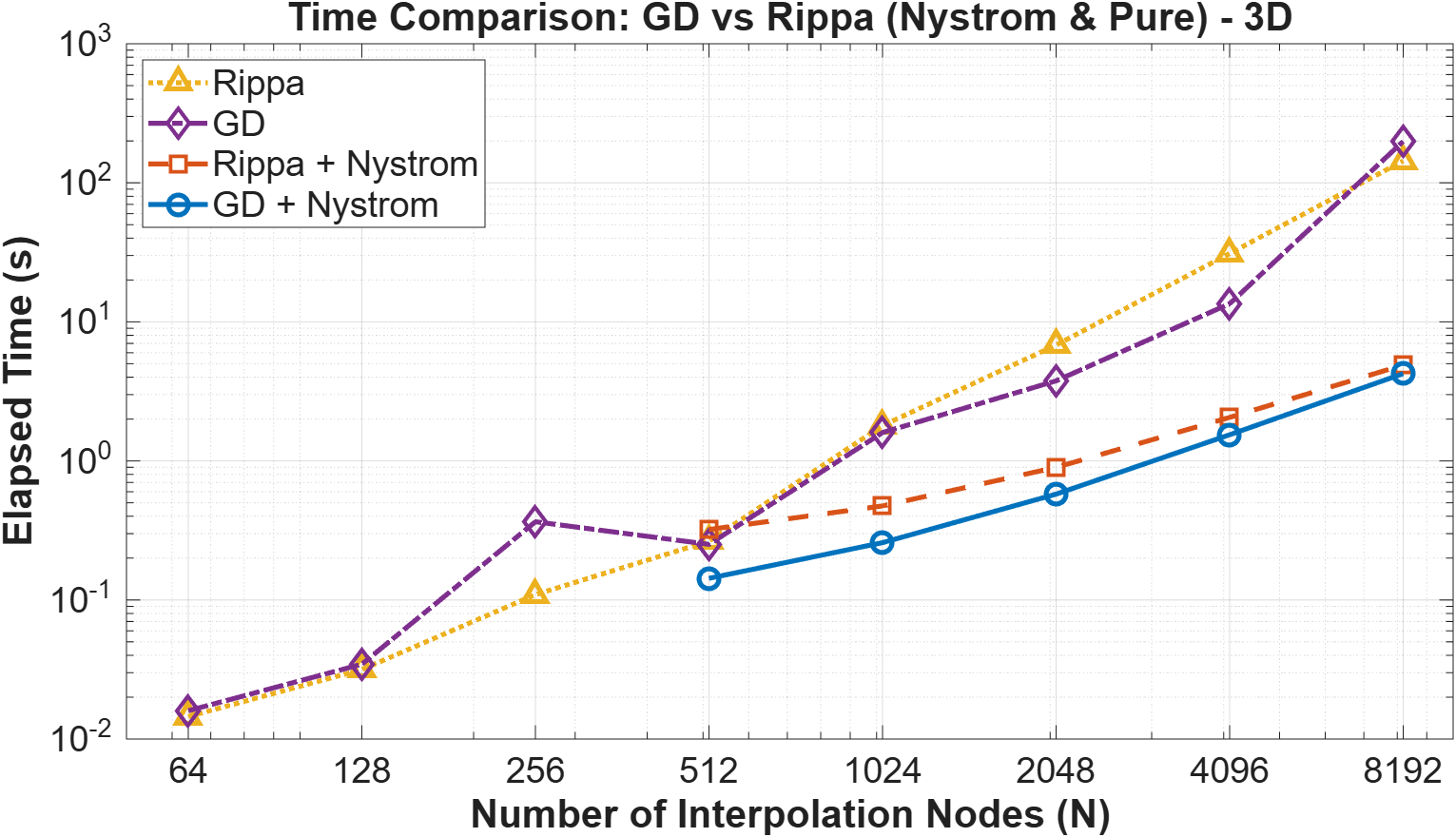}
        \caption{Time comparison}
        \label{fig:f7b}
    \end{subfigure}

    \caption{Performance plot for function $f_{6}$}
    \label{fig:f7plot}
\end{figure}

The test function $f_{8}: \Omega_{3} \rightarrow \R$ is \[f_{8} (\mathbf{x}) = \sin(2\pi (x_{1}^2 + 2x_{2}^2)) - \sin(2\pi(2x_{1}^2 +(x_{2}-0.5)^2 + x_{3}^2)),\]
$(x_{1},x_{2},x_{3}) \in [0,1] \times [0,1] \times [0,1]$.\\

\begin{figure}[H]
    \centering
    \begin{subfigure}[b]{0.7\textwidth} 
        \centering
        \includegraphics[width=\linewidth]{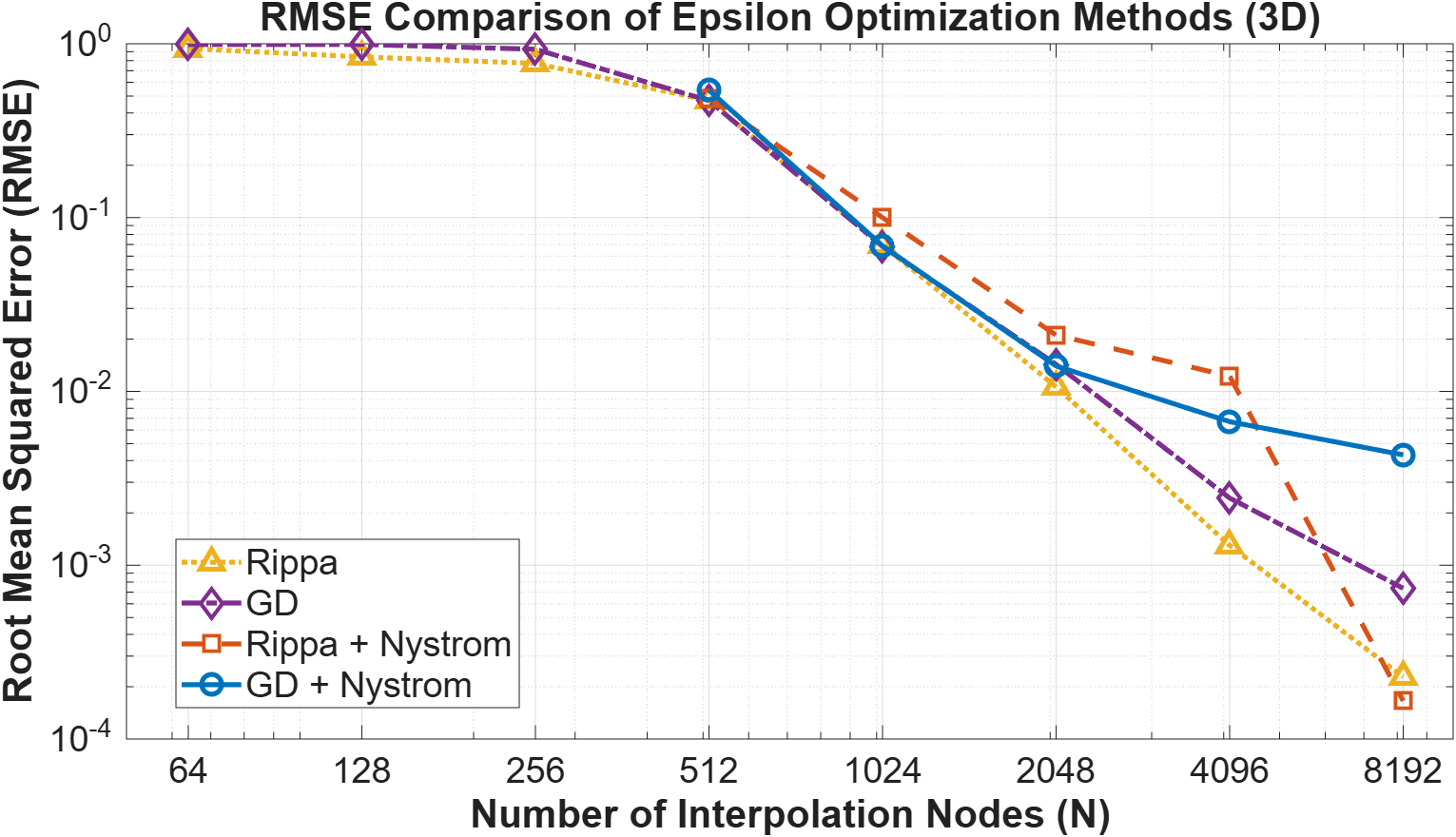}
        \caption{Accuracy Comparison}
        \label{fig:f8a}
    \end{subfigure}
    
    \vspace{1em} 

    \begin{subfigure}[b]{0.7\textwidth}
        \centering
        \includegraphics[width=\linewidth]{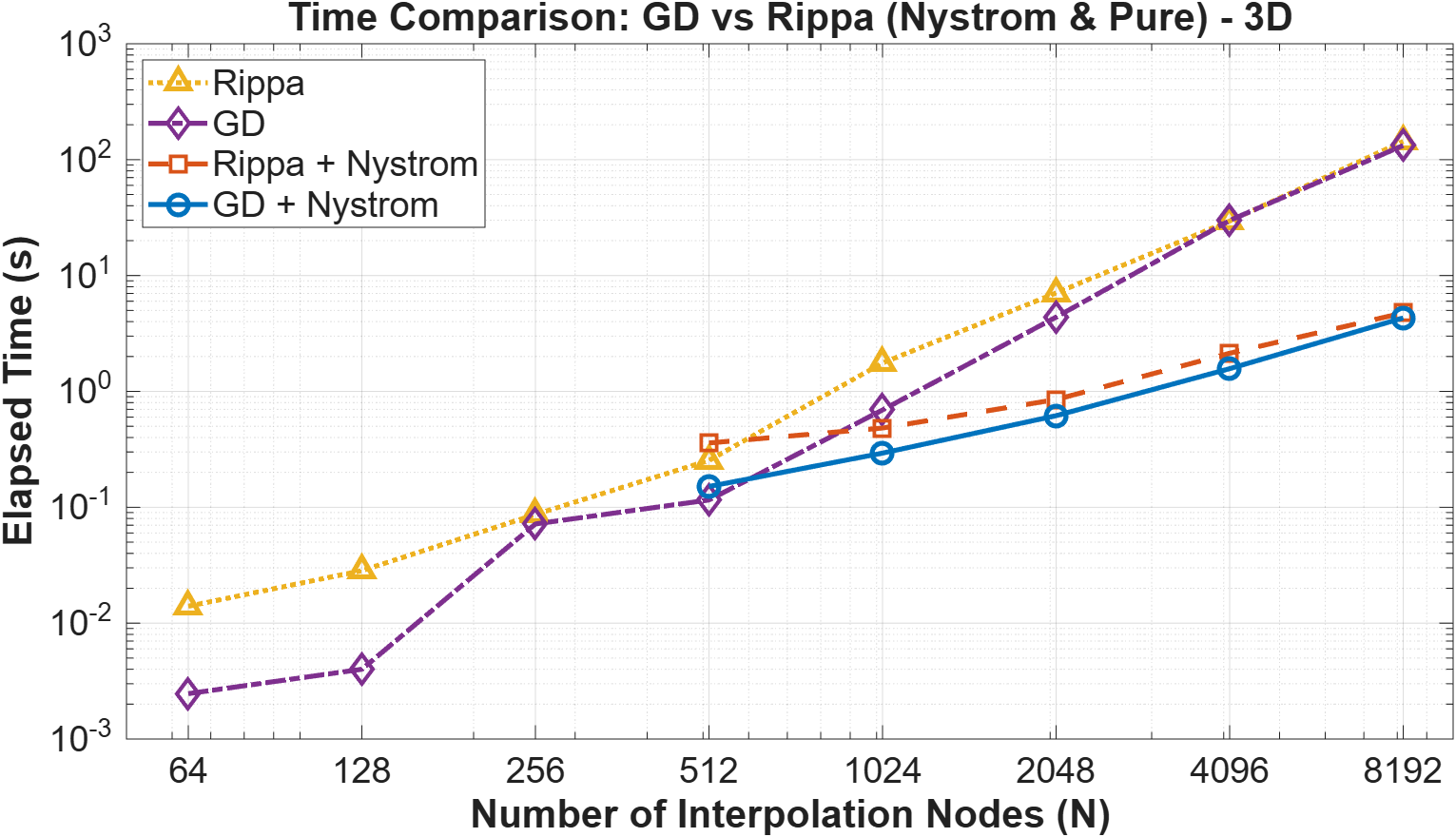}
        \caption{Time comparison}
        \label{fig:f8b}
    \end{subfigure}

    \caption{Performance plot for function $f_{8}$}
    \label{fig:f8plot}
\end{figure}

The timing comparisons (Figures \ref{fig:f6b}, \ref{fig:f7b}, \ref{fig:f8b}) strongly validate the theoretical complexity analysis. The accelerated methods (``GD + Nystr\"{o}m" and ``Rippa + Nystr\"{o}m") exhibit linear scaling with respect to the dataset size $N$. In particular, the gradient-based approach outperforms the grid-based search, reducing the runtime by efficiently navigating the parameter space with fewer objective evaluations.

In terms of accuracy (RMSE), the Nystr\"{o}m-based approximations track the behavior of the full LOOCV method closely. While the low-rank approximation ($m=200$) introduces a slight saturation effect at larger $N$ for the highly oscillatory function $f_8$, it maintains competitive accuracy for smoother functions. Crucially, the cosine mapping strategy successfully prevents the catastrophic divergence often seen in 3D RBF interpolation with uniform points, validating the robustness of our experimental setup.

\section{Conclusion}
\label{sec:conclusion}

In this work, we addressed the computational bottleneck of selecting the optimal shape parameter in RBF interpolation. By combining the Nystr\"{o}m low-rank approximation with the Woodbury matrix identity, we derived an efficient surrogate for the standard LOOCV objective function. This formulation successfully reduces the computational cost of a single evaluation from $\mathcal{O}(N^3)$ to $\mathcal{O}(Nm^2 + m^3)$, achieving linear scalability with respect to the dataset size $N$ once the landmark set size $m$ is fixed.

Our extensive numerical experiments across one, two, and three dimensions yield two primary conclusions regarding efficiency and robustness. The proposed Nystr\"{o}m-based formulation consistently delivers the expected speedups. In high-dimensional settings (3D), where standard interpolation matrices are notoriously ill-conditioned, the combination of our fast objective function with a Chebyshev-type node distribution (cosine mapping) proved essential for preventing divergence. This validates that LOOCV-based tuning can be made practical for large-scale datasets without sacrificing the qualitative behavior of the full method. While the Gradient Descent (GD) strategy offered the fastest execution times by minimizing objective evaluations, it exhibited sensitivity to the ill-conditioned landscapes inherent in RBF interpolation. In contrast, the grid-based ``Rippa with Nystr\"{o}m" approach provided a superior balance of speed and reliability. It retains the linear scalability of the approximation while avoiding the local minima and numerical instability that occasionally trapped the gradient-based solver.

These findings suggest that improving the optimization landscape, rather than merely reducing the cost of evaluating it, is essential for further progress. Future work must focus on developing more advanced optimization techniques, such as second-order methods, that are specifically designed to handle the ill-conditioned landscapes introduced by low-rank approximations. Additionally, exploring randomized sampling techniques, such as column-norm sampling \cite{mcdonaldpca} or CUR decomposition \cite{Wang2013}, could further reduce the $\mathcal{O}(NmdI_k)$ setup cost of k-means++ clustering. 

\appendix
 \section{NMI Computation Details}
\label{app:nmi}
This appendix provides a detailed description of the Normalized Mutual Information (NMI) 
computation procedure used in Section~2.2 to validate the stability of our landmark 
selection procedure.

\subsection{Motivation}
Since the k-means++ algorithm involves randomness in its initialization, one might be 
concerned that fixing a single random seed could inadvertently bias the numerical results. 
To verify that our conclusions are independent of initialization, we evaluate the stability 
by comparing multiple clustering results obtained with different random seeds.

\subsection{Procedure}
For each dataset size $N$, landmark count $K$, and dimension $D$, we perform the following steps:

\paragraph{Step 1: Generate Multiple Clusterings.} 
We run the k-means++ algorithm $R = 10$ times with different random initialization seeds 
(seeds $1, 2, \ldots, 10$). Each run produces a cluster assignment vector $C_r$ for 
$r = 1, 2, \ldots, R$, where $C_r \in \{1, 2, \ldots, K\}^N$ and $C_r(i)$ denotes the 
cluster assignment of the $i$-th data point.

\paragraph{Step 2: Compute Pairwise NMI.} 
We compute the NMI between every pair of clustering results. For clusterings $C_i$ and $C_j$ 
where $i < j$, the NMI is given by:
\[
NMI(C_i, C_j) = \frac{I(C_i, C_j)}{\sqrt{H(C_i)H(C_j)}},
\]
where $I(C_i, C_j)$ is the mutual information and $H(C_i)$, $H(C_j)$ are the entropies.
This yields $\binom{R}{2} = \binom{10}{2} = 45$ pairwise NMI values for each configuration.

The NMI takes values in $[0, 1]$ and measures the agreement between two clusterings:
\begin{itemize}
    \item $\text{NMI} = 1$: Perfect agreement (identical cluster assignments)
    \item $\text{NMI} = 0$: No agreement beyond random chance
    \item $\text{NMI} \approx 1$: High stability and reproducibility
\end{itemize}

\paragraph{Step 3: Aggregate Statistics.}
From the 45 pairwise NMI values, we compute:
\begin{itemize}
    \item \textbf{Mean NMI}: The arithmetic mean $\bar{\text{NMI}} = \frac{1}{45}\sum_{i<j} \text{NMI}(C_i, C_j)$
    \item \textbf{Standard deviation}: $\sigma_{\text{NMI}}$, measuring variability across pairs
\end{itemize}

\subsection{Interpretation}
A high mean NMI (e.g., $\bar{\text{NMI}} > 0.95$) with low standard deviation 
($\sigma_{\text{NMI}} < 0.05$) indicates that:
\begin{enumerate}
    \item The k-means++ algorithm produces consistent clusterings despite different initializations
    \item The landmark selection procedure is stable and reproducible
    \item Our numerical results do not depend on the particular choice of random seed
\end{enumerate}

By repeating this analysis across different dataset sizes $N \in \{1024, 2048, 4096, 8192\}$, 
landmark counts $K \in \{50, 100, 200, 300, 400\}$, and dimensions $D \in \{1, 2, 3\}$, 
we systematically demonstrate that stability improves as the ratio $N/K$ (data points per 
landmark) increases.

\subsection{Example}
Consider a concrete example with $K = 200$ landmarks and $N = 4096$ data points in $D = 2$ dimensions:
\begin{enumerate}
    \item Run k-means++ 10 times, producing clusterings $C_1, C_2, \ldots, C_{10}$
    \item Compute all pairwise NMIs: $\text{NMI}(C_1, C_2), \text{NMI}(C_1, C_3), \ldots, \text{NMI}(C_9, C_{10})$
    \item Obtain 45 NMI values, e.g., $\{0.978, 0.985, 0.981, \ldots\}$
    \item Calculate statistics: $\bar{\text{NMI}} = 0.982$, $\sigma_{\text{NMI}} = 0.008$
    \item Conclusion: Highly stable clustering with minimal sensitivity to initialization
\end{enumerate}

\section{K-means++ Algorithm for Landmark Selection}
\label{app:kmeans}

The k-means++ algorithm \cite{Arthur2007kmeans} is used to select $m$ landmark points from the interpolation node set $X = \{x_1, \ldots, x_N\}$. The algorithm proceeds as follows:

\begin{algorithm}[H]
\caption{K-means++ Landmark Selection}
\label{alg:kmeans_landmarks}
\begin{algorithmic}[1]
\REQUIRE Interpolation nodes $X = \{x_1, \ldots, x_N\} \subset \mathbb{R}^d$, number of landmarks $m$
\ENSURE Landmark indices $\mathcal{L} = \{\ell_1, \ldots, \ell_m\}$

\STATE \textbf{Initialization:}
\STATE Choose first cluster center $c_1$ uniformly at random from $X$
\STATE Set $\mathcal{C} = \{c_1\}$

\FOR{$k = 2$ to $m$}
    \FOR{each $x_i \in X$}
        \STATE $D(x_i) = \min_{c_j \in \mathcal{C}} \|x_i - c_j\|^2$
    \ENDFOR
    
    \STATE Choose $c_k \in X$ with probability $\frac{D(x_i)}{\sum_{j=1}^{N} D(x_j)}$
    \STATE $\mathcal{C} = \mathcal{C} \cup \{c_k\}$
\ENDFOR

\REPEAT
    \FOR{each $x_i \in X$}
        \STATE $\text{cluster}(x_i) = \arg\min_{j=1,\ldots,m} \|x_i - c_j\|$
    \ENDFOR
    
    \FOR{$j = 1$ to $m$}
        \STATE $c_j = \frac{1}{|\{i : \text{cluster}(x_i) = j\}|} \sum_{i: \text{cluster}(x_i) = j} x_i$
    \ENDFOR
\UNTIL{centers converge or maximum iterations reached}

\FOR{$j = 1$ to $m$}
    \STATE $\ell_j = \arg\min_{i=1,\ldots,N} \|c_j - x_i\|$
\ENDFOR

\STATE \textbf{return} $\mathcal{L} = \{\ell_1, \ldots, \ell_m\}$
\end{algorithmic}
\end{algorithm}

\subsection*{Key Properties}

\begin{enumerate}
    \item \textbf{Improved initialization:} K-means++ initialization (lines 2--11) provides provable approximation guarantees. Arthur and Vassilvitskii \cite{Arthur2007kmeans} showed that the expected cost is $O(\log m)$-competitive with the optimal clustering.
    
    \item \textbf{Distributed coverage:} By selecting centers proportional to squared distance from existing centers, k-means++ naturally spreads cluster centers across the data space, avoiding poor local minima common in random initialization.
    
    \item \textbf{Refinement:} The Lloyd's algorithm iterations (lines 12--22) further optimize the cluster assignments, though k-means++ initialization already provides a strong starting point.
    
    \item \textbf{Landmark mapping:} The final mapping step (lines 24--26) ensures that landmarks correspond to actual interpolation nodes rather than arbitrary cluster centroids, which is essential for constructing the Nyström approximation matrices $\mathbf{C}$ and $\mathbf{W}$.
\end{enumerate}

\subsection*{Implementation Details}

In our implementation, we use the following settings:
\begin{itemize}
    \item \textbf{Number of replicates:} 5 independent runs with different random initializations
    \item \textbf{Maximum Lloyd iterations:} 200
    \item \textbf{Convergence criterion:} Stop when center movement is below $10^{-6}$ or maximum iterations reached
    \item \textbf{Empty cluster handling:} If a cluster becomes empty during refinement, assign it to the point farthest from any existing center
\end{itemize}

The computational cost of k-means++ is $O(NmK\cdot d)$, where $K$ is the number of Lloyd iterations and $d$ is the spatial dimension. In practice, with $m \ll N$ and moderate $K$ (typically $< 100$), this cost is negligible compared to the savings achieved in subsequent LOOCV evaluations.

\bibliographystyle{siam}
\bibliography{mybib}
\end{document}